\documentstyle[11pt]{article}

\newcommand{\bplus}{\frame{$+$}}
\newcommand{\btimes}{\frame{$\times$}}
\newcommand{\freestar}{ \framebox[7pt]{$\star$} }
\newcommand{\bbox}{\Box}

\newcommand{\miab}{\mu_{i(ab-ba)}^{\bplus 1/2}}

\newcommand{\twtw}{ \widetilde{tw} }

\newcommand{\A}{ {\cal A} }

\newcommand{\C}{ \mbox{\bf C} }
\newcommand{\F}{{\cal F}}

\newcommand{\N}{ \mbox{\bf N} }
\newcommand{\R}{ \mbox{\bf R} }

\newcommand{\X}{ {\cal X} }

\newcommand{\Z}{ \mbox{\bf Z} }
\newcommand{\ee}{ \varepsilon }
\newcommand{\ff}{\varphi}

\newcommand{\lb}{\lbrack}
\newcommand{\rb}{\rbrack}

\newcommand{\aiqaik}{ a_{i_{1}} \cdots a_{i_{n}} }
\newcommand{\Xiqxin}{ X_{i_{1}} \cdots X_{i_{n}} }
\newcommand{\ziqzin}{ z_{i_{1}} \cdots z_{i_{n}} }
\newcommand{\iunun}{ i_{1} , \ldots ,i_{n} }
\newcommand{\Punun}{ P_{1} , \ldots ,P_{n} }
\newcommand{\Qunun}{ Q_{1} , \ldots ,Q_{n} }
\newcommand{\aunun}{ a_{1}, \ldots ,a_{m} }

\newcommand{\Bunur}{ B_{1}, \ldots ,B_{r} }
\newcommand{\Xunun}{ X_{1}, \ldots ,X_{m} }

\newcommand{\ncps}{ ( {\cal A} , \varphi ) }
\newcommand{\pstild}{ ( \widetilde{\cal A} , \widetilde{\varphi} ) }
\newcommand{\noncom}{ \mbox{non-commutative probability space} }

\newcommand{\ncpol}{ \C \langle X_{1} , \ldots ,X_{m} \rangle } 
\newcommand{\coefn}{ [ \mbox{coef } (n) ] }
\newcommand{\coefnpi}{ [ \mbox{coef } (n) ; \pi ] }
\newcommand{\coefnnn}{ [ \mbox{coef } ( \iunun ) ] }
\newcommand{\coefnnnpi}{ [ \mbox{coef } ( \iunun ) ; \pi ] }

\newcommand{\ecpi}{ \stackrel{\pi}{\sim} }
\newcommand{\ecrho}{ \stackrel{\rho}{\sim} }
\newcommand{\ecpiprim}{ \stackrel{\pi '}{\sim} }
\newcommand{\ecrhoprim}{ \stackrel{\rho '}{\sim} }
\newcommand{\ectrp}{ \stackrel{ \tau ( \rho ')}{\sim} }
\newcommand{\egdef}{:=}
\newcommand{\ecdef}{ \stackrel{def}{ \Leftrightarrow } }
\newcommand{\nn}{ \{ 1, \ldots ,n \} }

\newcommand{\Atild}{ \widetilde{\cal A} }
\newcommand{\atild}{ \widetilde{a} }
\newcommand{\btild}{ \widetilde{b} }
\newcommand{\comm}{ ab-ba }
\newcommand{\comtild}{ \atild \btild - \btild \atild }
\newcommand{\strings}{ \{ 1,2 \} ^{n} }
\newcommand{\lunun}{ l_{1} , \ldots ,l_{n} }
\newcommand{\infsum}{ \sum_{n=1}^{\infty} }
\newcommand{\rmua}{ R( \mu_{a} ) }
\newcommand{\mmub}{ M( \mu_{b} ) }
\newcommand{\ER}{R_E}
\newcommand{\er}{\ER}
\newcommand{\RE}{\ER}

\begin{document}

\title{\bf Commutators of free random variables}

\author{
\hbox{\centering
\vtop{\hsize 170pt\baselineskip=14pt
\begin{center}
Alexandru Nica\thanks{Partially supported by a Canada International 
Fellowship and by NSF Grant DMS-9622798}\\
Department of Mathematics\\
University of Michigan\\
Ann Arbor, MI 48109, U.S.A.\\
andu@math.lsa.umich.edu
\end{center}}
\
\vtop{\hsize 205pt\baselineskip=14pt
\begin{center}
Roland Speicher\thanks{Supported by a Heisenberg Fellowship of the DFG}\\
Institut f\"ur Angewandte Mathematik\\
Universit\"at Heidelberg\\
D-69120 Heidelberg, Germany \\
roland.speicher@urz.uni-heidelberg.de
\end{center}}
}}

\date{ }

\maketitle

\setlength{\baselineskip}{14pt}

\begin{abstract}
Let $\A$ be a unital $C^*$-algebra, given together with a
specified state $\ff:\A\to\C$. Consider two selfadjoint
elements $a,b$ of $\A$, which are free with respect to
$\ff$ (in the sense of the free probability theory of
Voiculescu). Let us denote $c:=i(ab-ba)$, where the $i$
in front of the commutator is introduced to make $c$
selfadjoint. In this paper we show how the spectral
distribution of $c$ can
be calculated from the spectral distributions of $a$ and
$b$. 
Some properties of the corresponding operation on
probability measures are also discussed.
The methods we use are combinatorial, based on the
description of freeness in terms of non-crossing partitions;
an important ingredient is the notion of $R$-diagonal pair,
introduced and studied in our previous paper \cite{NS3}.
\end{abstract}

$\ $

$\ $

\setcounter{section}{1}
\setcounter{equation}{0}
{\large\bf 1. Introduction and presentation of the results} 

$\ $

In this paper
we show how the combinatorial description of freeness can 
be used to determine
the distribution of
the commutator of two free random variables.

The concept of freeness was introduced by Voiculescu \cite{V1}
as a tool for studying free products of operator algebras.
It soon became clear that 
it is a promising point of view to consider freeness 
as a non-commutative analogue of the classical
probabilistic notion of independence, and this led to the development
of a free probability theory (see the monograph \cite{VDN},
or the recent survey in \cite{ICM}).

Let $a$ and $b$ be free random variables. Two basic 
problems in free probability theory, both solved 
by Voiculescu, consist in the study of $a+b$ and $ab$. The distributions
of these new random variables are given by the additive 
($\bplus$) and respectively multiplicative ($\btimes$)
free convolution of the distributions of $a$ and $b$. 
Voiculescu \cite{V2,V3} provided with the $R$- and $S$-transform two 
efficient analytical tools for dealing with $\bplus$ and $\btimes$.

To be more precise, let us recall that  if $(\A,\ff)$ is a 
{\it non-commutative probability space} -- i.e. $\A$ is a unital algebra, 
endowed with a linear functional $\ff:\A\to\C$ such that $\ff(1)=1$ --
then the {\it distribution} $\mu_a$ of an element $a\in\A$ is
the linear functional on $\C\lb X\rb$ defined by
$\mu_a(f)=\ff(f(a))$ for all $f\in\C\lb X\rb$. 
For every such distribution $\mu:\C\lb X\rb\to\C$ 
($\mu$ linear, $\mu(1)=1$)
one defines its $R$-transform $R(\mu)$ and its $S$-transform
$S(\mu)$ as special formal power series in an indeterminate $z$
in such a way that, for $a$ and $b$ free in
some $(\A,\ff)$, we have the formulas
\begin{equation}
R(\mu_{a+b})=R(\mu_a)+R(\mu_b) 
\end{equation}
and
\begin{equation}
S(\mu_{ab})=S(\mu_a)\cdot S(\mu_b),\qquad\mbox{if $\ff(a)\not=
0\not=\ff(b)$}.
\end{equation}
In the case of selfadjoint operators, where the distributions of
the random variables can be viewed as probability measures, the $R$-
and the $S$-transforms 
are related with the corresponding Cauchy transforms.
(For the definition of freeness, see Definition 2.1.4 below -- or,
for more details, the Section 2.5 of \cite{VDN}. For the definitions
of $R(\mu)$ and $S(\mu)$, as given by Voiculescu, we refer the reader 
to \cite{VDN}, Sections 3.2, 3.6.)

Having solved the problems of free addition and free multiplication,
the next \lq canonical' open problem is
the one of the free commutator, i.e. of describing the
distribution $\mu_{ab-ba}$ in terms of $\mu_a$ and $\mu_b$.
The goal of this paper is to present a solution to
the free commutator problem.

We shall use the alternative, combinatorial,  
approach to free random variables which was
developed by Speicher \cite{S1} and Nica and Speicher \cite{NS1,NS2,NS3}.
The main observation of \cite{S1} was that the coefficients of the
$R$-transform 
$R(\mu)$ -- which we call {\em free cumulants}
of $\mu$ -- can be
calculated via a precise combinatorial prescription using the lattice
of non-crossing partitions. Furthermore, in \cite{NS1,NS2}
it was also shown that the multiplicative free convolution can
be described directly in terms of the $R$-transform by
\begin{equation}
R(\mu_{ab})=R(\mu_a)\,\freestar\, R(\mu_b);
\end{equation}
here $\freestar$ is a combinatorial convolution of formal
power series, defined again by a precise combinatorial prescription
using the structure of non-crossing partitions.
The connection with the $S$-transform is achieved 
(see \cite{NS1}) by establishing
a combinatorial Fourier transform $\F$ which converts the 
$R$-transform into the $S$-transform,
\begin{equation}
S(\mu)=\F(R(\mu)),
\end{equation}
and which converts $\freestar$ into multiplication of power series,
\begin{equation}
\F(f\,\freestar\, g)=\F(f)\cdot \F(g).
\end{equation}

One of the most important advantages of this
combinatorial approach is that it can be
generalized in a straightforward way to multi-dimensional situations
\cite{S1,NS2}.
But this is exactly what is needed for the commutator problem, 
because
the crucial point turns out to consist in understanding, for $a$ 
and $b$ free,
the relation between $ab$ and
$ba$, i.e. 
in describing the two-dimensional joint distribution of
the pair $(ab,ba)$.

In the same spirit as in Eqns. (1.1), (1.3), we will approach
the distribution of the free commutator in terms of
its $R$-transform. 
The main result of the paper, Theorem 1.2, gives a description of
this $R$-transform in terms of our combinatorial convolution
$\freestar$; two other alternative descriptions, derived from the
one in Theorem 1.2, will be presented in the Corollaries 1.4 and
1.6.
We prefer to work with $i(ab-ba)$,
since, in a $C^*$-framework, 
this element is selfadjoint if $a$ and $b$ are so.

$\ $

{\bf 1.1 Notations.}
Let $\mu:\C\lb X\rb\to\C$ be a 
{\em distribution} (linear functional, normalized
by $\mu(1)=1$) and let 
$\lb R(\mu)\rb(z)=\sum_{n=1}^\infty \alpha_nz^n$
be its $R$-transform.
We shall denote by $\ER(\mu)$ the generating series of
the even free cumulants $(\alpha_{2n})_{n=1}^\infty$, i.e.
\begin{equation}
\lb \ER(\mu)\rb(z):=\sum_{n=1}^\infty \alpha_{2n}z^n.
\end{equation}

$\ $

In general, $\ER(\mu)$ contains less information than $R(\mu)$;
however, in the case when $\mu$ is {\em even} -- i.e. 
$\mu(X^{2n+1})=0$ for all $n$ and hence  
also $\alpha_{2n+1}=0$ for all $n$ -- then $R(\mu)$ can
be recaptured from $\ER(\mu)$ by
\begin{equation}
\lb R(\mu)\rb(z)=\lb\ER(\mu)\rb(z^2).
\end{equation}

$\ $

{\bf 1.2 Theorem.}
Let $\ncps$ be a $\noncom$, let $a,b$ be in $\A$ and
consider the
element $c:=i(ab-ba)$.
If $a$ is free from $b$ then
\newline
i) $\mu_c$ is even, and\newline
ii) we have the equation
\begin{equation}
\ER(\mu_c)=2(\ER(\mu_a)\,\freestar\,\ER(\mu_b)\,\freestar\, Zeta),
\end{equation}
where $Zeta$ is the special series
\[
Zeta(z)\egdef \sum_{n=1}^\infty z^n.
\]

$\ $

Of course, parts i) 
and ii) of our theorem can be combined into only one formula:
\begin{equation}
\lb R(\mu_{i(ab-ba)})\rb(z)=2(\ER(\mu_a)
\thinspace\freestar\thinspace \ER(\mu_b)\thinspace\freestar 
\thinspace Zeta)(z^2).
\end{equation}

Note that Theorem 1.2 does not only state that the 
distribution of the commutator is always even, 
but also, even more strikingly, that the distribution of the commutator
depends only on the even 
free cumulants of $a$ and on the even free cumulants
of $b$.

It is not mandatory to keep the result of Theorem 1.2 
formulated in terms of $\freestar$, since (modulo 
a technical detail) we can apply in (1.8) the combinatorial
Fourier transform and turn $\freestar$ into usual multiplication.
Before doing this let us first present a direct application of
formula (1.8) itself.

$\ $

{\bf 1.3 Application (commutator with symmetric Bernoulli distribution).}
Let $a$ and $b$ be free in some non-commutative probability space
$(\A,\ff)$, such that $b$ has distribution
$\mu_b=1/2(\delta_{-1}+\delta_{+1})$. 
Then we have:
\begin{equation}
\mu_{i(ab-ba)}= \mu_a\,\bplus \,\mu_{-a}.
\end{equation}

$\ $

Indeed, one finds out in this case (for instance by plugging $\mu_{b}$
into the Eqn. (5.2) of Proposition 5.2 below) that $\er ( \mu_{b} )$ is 
equal to $Moeb,$ the inverse of $Zeta$ under $\freestar$. Therefore
(1.9) reduces to
\[
[ R( \mu_{i(ab-ba)} ) ] (z) \ = \
2 [ \er ( \mu_{a} ) ] ( z^{2} ) \ = \ 
[ R( \mu_{a} ) ] (z) + [ R( \mu_{-a} ) ] (z) \ = \
[ R( \mu_{a} \, \bplus \, \mu_{-a} ) ] (z),
\]
which implies (1.10).

$ \ $

We now give the `analytical' reformulation of Theorem 1.2, which was 
announced prior to 1.3. Recall that if $\mu : \C [X] \rightarrow \C$
is a linear functional with $\mu (1) = 1,$ then the series
\[
\lb M(\mu)\rb(z) \egdef \sum_{n=1}^\infty \mu(X^n)z^n
\]
is called the {\it moment series} of $\mu$. By applying the combinatorial
Fourier transform in (1.8) and by taking into account the relations which
exist between the series $R(\mu)$, $S(\mu)$, $M(\mu)$, one obtains
the following.

$\ $

{\bf 1.4 Corollary.}
Let $\ncps$ be a $\noncom$ and let $a,b\in\A$ be  free and with non-zero
variances ($\ff(a^2)-(\ff(a))^2 \neq 0 \neq \ff(b^2)-(\ff(b))^2$). 
Then the moment series of $i(ab-ba)$ is given by
\begin{equation}
\lb M(\mu_{i(ab-ba)})\rb(z)=
\Bigl( \frac{2}{ w(1+ \frac w2 )(1+w)^{2} } \cdot 
\lb \er ( \mu_a ) \rb^{<-1>} ( \frac w2 ) \cdot 
\lb \er ( \mu_b ) \rb^{<-1>} ( \frac w2 ) 
 \Bigr)^{<-1>} ( z^2 ),
\end{equation}
where the notation `$h^{<-1>}$' is used for the inverse under composition
of a formal power series $h$ without constant term, and 
with non-zero linear
term $(h(0) = 0 \neq h' (0)).$

$\ $

This corollary can be thought of as an algorithm which tells how to 
obtain the moment series of $i(ab-ba)$ from the ones of $a$ and $b$, 
via a succession of algebraic operations, simple substitutions, and 
inversions under composition. Strictly speaking, the Equation (1.11)
only connects $M( \mu_{i(ab-ba)} )$ to the sequences of even free 
cumulants of $a$ and $b$, and not to $M(\mu_a)$ and $M(\mu_b)$; but 
it is known that the free cumulants are obtained from the moments by 
using the same types of operations -- see, e.g., \cite{VDN}, Section 3.3.

Let us now examine some concrete examples of how Corollary 1.4 can
be applied. The first step of a calculation based on 1.4 is to find out
the inverse $\er$-transforms $[ \er ( \mu_{a} ) ]^{<-1>}$ and
$[ \er ( \mu_{b} ) ]^{<-1>}$. The inverse $\er$-transforms of a few 
distributions which occur frequently in free probability are listed in
the following table.

\begin{tabular}{|c|c|}     \hline
                              &                                   \\
Distribution $\mu$            &  $[ R_{E} ( \mu ) ]^{<-1>} (w)$   \\ 
                              &                                   \\  \hline
                              &                                   \\    
semicircular of radius $r>0$  &                                   \\
                              &                                   \\
$\left(  \begin{array}{c}
{d \mu (t)=2 \pi^{-1} r^{-2} (r^{2} -t^{2})^{1/2} dt \mbox{ on } [-r,r]} \\
{[ R( \mu ) ] (z) = r^{2} z^{2} /4  \ - \ \mbox{see \cite{VDN}, 3.4.2}}
\end{array}   \right)$        &    $\begin{array}{c}  4w  \\
                                   \overline{r^{2}} \end{array}$  \\
                              &                                   \\  \hline  
                              &                                   \\    
free Poisson of parameters $\alpha , \beta > 0$  
                              &                                   \\
$\left( [ R( \mu ) ] (z) = \alpha \beta z / (1- \beta z) \ - \  
\mbox{see \cite{VDN}, 3.7} 
\right)$                      &  $\begin{array}{c} w \\ 
                                 \overline{ \beta^{2} 
                                 ( \alpha +w)} \end{array}$        \\
                              &                                    \\  \hline  
                              &                                    \\    
arcsine law on $[-r,r]$, $r>0$ &                                   \\
                              &                                    \\
$\left(  \begin{array}{c}
{d \mu (t) = \pi^{-1} (r^{2} - t^{2})^{-1/2} dt \mbox{ on } [-r,r]} \\
{[ R( \mu )](z)=-1+ \sqrt{1+r^{2}z^{2}}
 \  - \  \mbox{see \cite{VDN}, 3.4.5}}
\end{array}   \right)$        &   $\begin{array}{c} 
                                  \underline{2w+w^{2}} \\
                                  r^{2}  \end{array}$               \\
                              &                                     \\  \hline  
                              &                                     \\    
Bernoulli                     &                                     \\
                              &                                     \\
$\left(  \begin{array}{c}
\mu =  \lambda \delta_{t_{0}} + (1- \lambda ) \delta_{t_{1}}   \\ 
\mbox{for some } t_{0} < t_{1} \mbox{ in } \R , \lambda \in (0,1)
\end{array}   \right)$        &  $\begin{array}{c} 
                                 w(1+w)(1+2w)^{2} \\
                                 \overline{(t_{1}-t_{0})^{2} (w^{2}
                                 +w+( \lambda - \lambda^{2} )) }           
                                 \end{array}$                        \\
                              &                                      \\ \hline
\end{tabular}

\[
\mbox{Table 1}
\]

By using Eqn. (1.11) one can in principle obtain the distribution of the
commutator of two free elements $a,b$ whenever $\mu_{a}$ and $\mu_{b}$ are
taken from the Table 1. A detail which cannot be ignored, though, is that 
the calculation of $M( \mu_{i(ab-ba)} )$ requires one more inversion under
composition; if we just restrict our attention to the situations in Table 1,
then, due to the fact that all the $[ \er ( \mu ) ]^{<-1>}$'s in the table
are rational functions, the remaining inversion under composition comes to 
solving a (possibly not too pleasant) algebraic equation.

Let us point out that, even though the Corollary 1.4 is formulated in a 
combinatorial framework, its applications take place in a $C^*$-context. 
If the non-commutative probability space $(\A,\ff)$ is actually
a {\em $C^*$-probability space} -- i.e., $\A$ is a unital 
$C^*$-algebra, and $\ff$ is a state on $\A$ -- then the
distribution $\mu_c$ of a selfadjoint element $c\in\A$ is
actually a probability measure with compact support on $\R$.
In such a case, the moment series $M(\mu_c)$ can still be viewed
as a formal power series, but it is also an analytic function
on a neighborhood of zero, and is related to the Cauchy transform
$G$ of $\mu_c$ via the formula
\begin{equation}
1+\lb M(\mu_c)\rb(z) \ = \ z^{-1} G(z^{-1}),\qquad
0<\vert z\vert <\Vert c\Vert^{-1}.
\end{equation}
The algebraic equation for $M( \mu_{i(ab-ba)} )$ obtained by using 
(1.11) can thus be converted (with the help of (1.12), written for
$c = i(ab-ba)$), into an algebraic equation satisfied by $G.$
From this point one can, in nice cases, solve for $G$ and then
recover the measure $\mu_{i(ab-ba)}$ itself by using the
Stieltjes inversion formula (see e.g. \cite{Akh}). Some situations of this 
kind -- when all the steps work out and yield $\mu_{i(ab-ba)}$ 
explicitly -- are shown next. 
We should also point out that in these examples the spectrum of $i(ab-ba)$
is in particular obtained, as the support of $\mu_{i(ab-ba)}$.

$\ $

{\bf 1.5 Examples.} In all the examples discussed here, $a$ and $b$ are
free selfadjoint elements in a $C^{*}$-probability space $\ncps$. The
element $i(ab-ba) \in \A$ is denoted by $c$, and the Cauchy transform of
$\mu_{c}$ is denoted by $G.$

\vspace{6pt}

1) Assume that $a$ is semicircular (say, for definiteness, that it has
radius 2) and that $b$ is a projection with $\varphi (b) =: \lambda \in
(0,1).$ Then (1.11) gives a quadratic equation for $M( \mu_{c} );$ when
converted into an equation for $G,$ this becomes:
\begin{equation}
\zeta^{2} G( \zeta )^{2} - 2 \zeta^{3} G( \zeta ) + (2 \zeta^{2} - 1
+ 4 \lambda (1- \lambda ) ) \ = \ 0, \ \ \  \zeta \in \C , \ 
\mbox{Im} \zeta > 0.
\end{equation}
By solving for $G$ and then by applying the Stieltjes inversion formula, 
we obtain:
\begin{equation}
\mu_{c} \ = \ \sqrt{1 - 4 \lambda (1- \lambda ) } \  \delta_{0} \ + \
\frac{1}{\pi |t|}
\sqrt{4 \lambda (1- \lambda ) - (t-1)^{2} } \ 
\chi_{ [ - \beta , - \alpha ] \cup [ \alpha , \beta ] }\, dt,
\end{equation}
where $\alpha := \sqrt{1- 2 \sqrt{\lambda (1- \lambda )} } \in [0,1)$
and $\beta := \sqrt{1+ 2 \sqrt{\lambda (1- \lambda )} } \in (1, 
\sqrt{2} ].$ Note that the spectrum of $c$ (which is 
$[ - \beta , - \alpha ] \cup \{ 0 \} \cup [ \alpha , \beta ]$) is an 
interval only if $\lambda = 1/2;$ in this case $\mu_{c}$ is semicircular
of radius $\sqrt{2}$, as one could also infer directly from Application 1.3.

\vspace{6pt}

2) Assume that both $a$ and $b$ are semicirculars of radius 2. Then (1.11)
gives us a cubic equation for $M( \mu_{c} )$; hence we also get via
(1.12) a cubic equation for $G$, which turns out to be:
\begin{equation}
\zeta G( \zeta )^{3} + G( \zeta )^{2} -  \zeta G( \zeta ) + 1 \ = \ 0.
\end{equation}
By solving for $G$ and then by using the Stieltjes inversion formula 
we obtain in this case that $\mu_{c}$ is absolutely continuous with respect
to Lebesgue measure. The support of $\mu_{c}$ is the interval 
$[ -r,r ]$, with $r = \sqrt{ (11+5 \sqrt{5} )/2 };$ its density is
\begin{equation}
\frac{d\mu (t)}{dt} \ = \ \frac{ \sqrt{3} }{2 \pi |t| } \left(
\frac{ 3t^{2} + 1}{ 9h(t) } - h(t) \right), \ \ \ |t| \leq
\sqrt{ (11+5 \sqrt{5} )/2 },
\end{equation}
where
\begin{equation}
h(t) \ = \
\sqrt[3]{ \frac{18t^{2} +1}{27} \ + \ 
\sqrt{\frac{t^{2}(1+11t^{2}-t^{4})}{27} } } .
\end{equation}

This example has an alternative direct derivation, as follows. One notes 
first that $c$ has the same distribution as the anti-commutator $ab+ba$
(see Proposition 1.10 below). Then one writes 
\begin{equation}
ab+ba \ = \
\left( \frac{a+b}{\sqrt{2}} \right) ^{2} \ - \
\left( \frac{a-b}{\sqrt{2}} \right) ^{2} ,
\end{equation}
and uses the fact that $(a+b)/ \sqrt{2}$ and $(a-b)/ \sqrt{2}$ are also
free semicircular of radius 2 ( see e.g. \cite{VDN}, Section 2.6). Since 
the square of a semicircular of radius 2 is Poisson of parameters 1,1
(with $R$-transform $Zeta (z) = z/(1-z),$ as listed in Table 1),
one gets:
\[
[ R( \mu_{c} ) ] (z) \ = \
[ R( \mu_{ab+ba} ) ] (z) 
\]
\begin{equation}
= \ [ R( \mu_{(a+b)^{2} /2} ) ] (z) + [ R( \mu_{(a-b)^{2} /2} ) ] (-z) 
\ = \ \frac{z}{1-z} + \frac{-z}{1+z} \ = \
\frac{2z^{2}}{1-z^{2}}.
\end{equation}
The formulas (1.15-17) then follow from (1.19) by using the direct relation
which exists between the $R$-transform and the Cauchy transform (see
\cite{VDN}, Section 3.3).

In view of (1.19), the distribution $\mu_{c}$ of this example could be
called `symmetric Poisson'. We note that (1.19) is also obvious from the
statement of Theorem 1.2. Indeed, in this case we have
$[ R( \mu_{a} ) ] (z) = [ R( \mu_{b} ) ] (z) = z$, which is the unit for
the operation $\freestar$; thus Eqn. (1.9) becomes
\[
[ R( \mu_{c} ) ] (z) \ = \ 2 \ Zeta ( z^{2} ) \ = \ 
\frac{2z^{2}}{1-z^{2}}.
\]

The nice simple trick of (1.18) does not seem to work in other examples,
because $(a+b)/ \sqrt{2}$ and $(a-b)/ \sqrt{2}$ are not free in general
(actually, the freeness of these elements implies that $a$ and $b$ are
semicircular -- see \cite{N}, Section 5).

\vspace{6pt}

3) Assume that both $a$ and $b$ are projections. Then the equations 
obtained for $M( \mu_{c} )$ and $G$ are of degree four, but can be solved
without much difficulty, because they are bi-quadratic. We only present
the formula for $\mu_{c}$ in two particular cases where the calculations
are nicer than in the generic one. For both the particular cases it is
convenient to denote
\begin{equation}
\xi \ := \ 4 \varphi (a) ( 1 - \varphi (a) ) \in (0,1]
\end{equation}
(where it is assumed that $\varphi (a) \neq 0,1).$

\vspace{4pt}

a) If $\varphi (b) = 1/2$, then (with $\varphi (a) \in (0,1)$ arbitrary
and $\xi$ as in (1.20)) we have:
\begin{equation}
\mu_{c} \ = \ \sqrt{1- \xi} \  \delta_{0} \ + \ 
\frac{1}{\pi |t|}
\sqrt{ \frac{ 4t^{2} - 1+ \xi }{ 1 - 4t^{2} } } \ 
\chi_{ [ - 1/2 , - \sqrt{1- \xi}/2  ] \cup
[ \sqrt{1- \xi}/2 , 1/2 ] }  dt.
\end{equation}
Here too, the spectrum of $c$ is an interval only if $\varphi (a) = 1/2;$
in this case $\mu_{c}$ is the arcsine law on $[ - 1/2 , 1/2 ].$

\vspace{4pt}

b) If $\varphi (a) = \varphi (b) =: \lambda \in (0,1),$ then there are
several subcases to consider. For $0 < \lambda \leq \frac{1}{2} -
\frac{1}{ \sqrt{8} }$ we have (with $\xi$ from (1.20)) that:
\begin{equation}
\mu_{c} \ = \ \sqrt{1- \xi} \  \delta_{0} \ + \ 
\frac{2}{\pi}
\sqrt{ \frac{ \xi (1- \xi ) - t^{2} }{ (1-4t^{2})
(1+ \sqrt{1-4t^{2}} )(1-2 \xi + \sqrt{1-4t^{2}} ) } } \
\chi_{ [ - \sqrt{\xi (1- \xi )} , \sqrt{\xi (1- \xi )} ] }  dt
\end{equation}
(which means in particular that the norm of $c$ decreases with $\lambda$).
For $\frac{1}{2} - \frac{1}{ \sqrt{8} } \leq \lambda \leq \frac{1}{2}$ 
we have:
\begin{equation}
\mu_{c} \ = \ \sqrt{1 - \xi} \ \delta_{0} \ + \ 
h_{1} (t) \chi_{ [ - \sqrt{\xi (1- \xi )} , \sqrt{\xi (1- \xi )} ] } \ dt 
\ + \ h_{2} (t) \chi_{ [ -1/2,  - \sqrt{\xi (1- \xi )} ] \cup
[ \sqrt{\xi (1- \xi )} , 1/2 ] } \ dt ,
\end{equation}  
where the densities $h_{1}, h_{2}$ are given by:
\[
h_{1} (t)  \  =  \  
\frac{1}{\pi}  \sqrt{ \frac{ 2 \xi -1+ \sqrt{1-4t^{2}} }{(1+
\sqrt{1-4t^{2}} )(1-4t^{2} ) } }                                
\]
\begin{equation}
\end{equation}
\[
h_{2} (t)  \  =  \  
\frac{1}{\pi}  \sqrt{ \frac{ 2t^{2} -1+ \xi + 2t 
\sqrt{t^{2} - \xi (1- \xi ) } }{ t^{2} (1-4t^{2} ) } } \  .                   
\]
The situations when $1/2 \leq \lambda < 1$ are obtained from those when
$0 < \lambda \leq 1/2$, by replacing $a$ and $b$ with $1-a$ and $1-b$,
respectively.

$ \ $

An element $a$ in a non-commutative probability space $\ncps$ is called
{\em even} if its distribution $\mu_{a}$ is so, i.e. if 
$\varphi ( a^{n} ) = 0$ for $n$ odd. In the particular case when we look
at the free commutator of two even elements, we have a third reformulation
of Eqn. (1.8) in Theorem 1.2, in terms of $S$-transforms.

$\ $

{\bf 1.6 Corollary.} Let $\ncps$ be a non-commutative probability space,
and let $a,b \in \A$ be even, free, and with non-zero variances 
$\gamma_{a} = \varphi ( a^{2} ) - ( \varphi (a))^{2} ,$
$\gamma_{b} = \varphi ( b^{2} ) - ( \varphi (b))^{2} .$ Consider the 
element $c := i(ab-ba) \in \A$, and its variance
$\gamma_{c} = \varphi ( c^{2} ) - ( \varphi (c))^{2} .$ Then
$\gamma_{c} = 2 \gamma_{a} \gamma_{b} ( \neq 0),$ and we have the equation
in $S$-transforms:
\begin{equation}
[ S( \mu_{ c^{2} / \gamma_{c} } ) ] (w) \ = \
\frac{ 1+ \frac w2 }{ 1+w } \cdot
[ S( \mu_{ a^{2} / \gamma_{a} } ) ] ( \frac w2 ) \cdot
[ S( \mu_{ b^{2} / \gamma_{b} } ) ] ( \frac w2 ).
\end{equation}

$ \ $

Note that the $S$-transform $S( \mu_{ c^{2} / \gamma_{c} } )$ determines
the distribution of $c^{2}$, which in turn determines the one of $c$
(because $c$ is even, by Theorem 1.2.i). Hence Equation (1.25) can also 
be viewed as a solution to the free commutator problem, in the 
particular situation when we work with even elements. The proof of
Corollary 1.6 goes in parallel with the one of 1.4 -- 
see Section 5.5 below.

$\ $

{\bf 1.7 Application (iterated free commutators)} Let $\ncps$ be a 
non-commutative probability space, and let $( a_{m} )_{m=1}^{\infty}$ be
a sequence of free elements of $\A$, all having the same even distribution
$\mu$. We form the sequence of iterated commutators:
\[
c_1 = a_1, \ \ \ c_m = i(c_{m-1} a_m - a_m c_{m-1} ), \ \ \ m \geq 2,
\]
and we consider the problem whether this sequence has a limit distribution
for $m \rightarrow \infty$.

Let $\gamma$ be the variance of $\mu$. From the formula for variances in
Corollary 1.6 it follows that the variance of $c_m$ is
$\frac{1}{2} (2 \gamma )^{m}$, $m \geq 1.$ Hence, in order for the limit
distribution of the $c_m$'s to exist and be non-trivial, we need to make 
the additional assumption that $\gamma = 1/2$ (which implies that $c_m$
has variance $1/2$ for all $m$).

The distribution $\mu_{ 2a_{m}^{2} }$ does not depend on $m$ (but only on
$\mu$, the common distribution of the $a_m$'s); we denote 
$S( \mu_{ 2a_{m}^{2} } ) =: g.$ Then the repeated application of 
Corollary 1.6 leads to the formula:
\begin{equation}
[ S( \mu_{ 2c_{m}^{2} } ) ] (w) \ = \ 
\Bigl[ \prod_{k=1}^{m-1} g( \frac{w}{2^k} ) \Bigr] \cdot
g( \frac{w}{ 2^{m-1} }) \cdot \frac{ 1+ \frac{w}{ 2^{m-1} } }{ 1+w },
\ \ \ m \geq 2.
\end{equation}
It is easy to check that the series $\prod_{k=1}^{m-1} g( \frac{w}{2^k} )$
have, for $m \rightarrow \infty$, a coefficient-wise limit distribution, 
which we will denote $\prod_{k=1}^{\infty} g( \frac{w}{2^k} )$. From (1.26)
it follows that $S( \mu_{ 2c_{m}^{2} } )$ converges coefficient-wise, for
$m \rightarrow \infty$, to
$\frac{1}{1+w} \cdot \prod_{k=1}^{\infty} g( \frac{w}{2^k} ).$ Since the
moments of a functional are recaptured from the coefficients of its 
$S$-transform via polynomial equations (see e.g. \cite{V3}), we can thus
infer the moment-wise convergence of the functionals
$ ( \mu_{ 2c_{m}^{2} } )_{m=1}^{\infty} $. Finally, due to the fact that
all the $c_m$'s are even, the convergence of the sequence 
$( \mu_{ 2c_{m}^{2} } )_{m=1}^{\infty} $ is equivalent to the one of the
sequence $( \mu_{ c_{m} } )_{m=1}^{\infty}.$ 

The conclusion is hence that: under the additional assumption that $\mu$
has variance $1/2$, the sequence of distributions 
$( \mu_{ c_{m} } )_{m=1}^{\infty}$ converges moment-wise, and moreover,
if $c_{\infty}$ is an element (in some non-commutative probability space)
such that $\mu_{c_{\infty}}$ = ${\mbox{lim}}_{m \rightarrow \infty}
\mu_{c_{m}}$, then we have:
\begin{equation}
[ S( \mu_{ 2c_{\infty}^{2} } ) ] (w) \ = \
\frac{1}{1+w} \cdot \prod_{k=1}^{\infty} g( \frac{w}{2^k} ). 
\end{equation}
From (1.27) it is clear that the limit distribution $\mu_{\infty}$ of the
$\mu_{c_{m}}$'s effectively depends on the input distribution $\mu$ we start
with. Some examples that can be easily calculated are:

\vspace{4pt}

a) If $\mu$ is the symmetric Bernoulli distribution
$\frac{1}{2} ( \delta_{-1/2} + \delta_{1/2} ),$ then $\mu_{\infty}$ is 
semicircular of radius $\sqrt{2}.$

\vspace{4pt}

b) If $\mu$ is the arcsine law $2 \pi^{-1} (1-4t^{2})^{-1/2}$ on
$[ -1/2 , 1/2 ]$, then $\mu_{\infty}$ is the symmetric (free)
Poisson distribution of Example 1.5.2, properly normalized with a 
dilation by $1/4$. This example, b), can be in fact viewed as a 
consequence of the previous one, a) -- see Section 1.17 below.

\vspace{4pt}

c) If $\mu$ is semicircular of radius $\sqrt{2}$, then the common
distribution of the elements $2a_{m}^{2}$ is free Poisson of parameters
1,1, and this implies that the function $g$ defined prior to (1.26) is
just $g(w) = 1/(1+w).$ The right-hand side of (1.27) becomes hence
$\prod_{k=0}^{\infty} 1/(1+ \frac{w}{2^{k}} );$ one can transform this
infinite product into a sum (see e.g. \cite{And}, Theorem 2.1), and 
rewrite it as ${\mbox{exp}}_{1/2} ( -2w ),$ where the 
$\frac{1}{2}$-exponential series is
\[
{\mbox{exp}}_{1/2} ( t ) \ = \ 1 + \sum_{n=1}^{\infty}
\frac{ t^{n} }{ [1]_{1/2} [2]_{1/2} \cdots [n]_{1/2} } ,
\]
with $[n]_{1/2} = 1 + 1/2 + \cdots + 1/2^{n-1} = 2- 1/2^{n-1}.$ It is not
clear though how $\mu_{\infty}$ explicitly looks in this case.

$\ $

Up to now, we have restricted our attention 
to the level of random variables.
But as it follows directly from general freeness considerations,
taking the commutator of free random variables is really an
operation on the level of distributions. 
The main result of the paper, Theorem 1.2, can be viewed as
a description of this operation in terms of the combinatorial
convolution $\freestar$ and the $\ER$-transform.
In the rest of the Introduction we will switch to this
alternative point of view and we
will examine the structure of this operation more closely.

$\ $

{\bf 1.8 Notation.} Let $\nu_1$, $\nu_2$ be two distributions (in the
algebraic sense  of Notations 1.1). We shall denote the free commutator 
of $\nu_1$ and $\nu_2$ by $\lb \nu_1 \Box \nu_2 \rb$, and their free 
anti-commutator by $\{ \nu_1 \Box \nu_2 \}$. More explicitly, if $a$ 
and $b$ are free random variables in some $(\A,\ff)$ such that 
$\nu_1=\mu_a$ and $\nu_2=\mu_b$, then
\begin{equation}
\lb \nu_1 \Box \nu_2 \rb \egdef \mu_{i(ab-ba)} , \ \ \ 
\mbox{ and } \ \ \ 
\{ \nu_1 \Box \nu_2 \} \egdef  \mu_{ab+ba}.
\end{equation}

$\ $

{\bf 1.9 Remark.} In the $C^*$-context, we thus have the operations 
of taking the free commutator and of taking the free anti-commutator 
for probability measures with compact support. In \cite{BV} it is
shown that one can extend these operations (or, more generally, 
operations given by non-commutative selfadjoint polynomials) from 
probability measures with compact support to all probability measures.

$\ $

In general, the commutator and the anti-commutator are of course
different. But one of the byproducts of our results of Section 4
is the fact that they coincide for even distributions.

$\ $

{\bf 1.10 Proposition.}
Let $\nu_1$ and $\nu_2$ be even distributions. Then the commutator
and the anti-commutator of $\nu_1$ and $\nu_2$ are the same,
\[
\lb \nu_1\Box\nu_2\rb= \{ \nu_1\Box\nu_2\}.  
\]

$\ $

The description of the free anti-commutator for general distributions
is still an open problem. Our general result for the free commutator
relies on the fact that due to cancellations we can reduce the
general case to the even case.

$\ $

Clearly, the operation of taking the free commutator is not
associative. But, Theorem 1.2 shows that only the factor 2 in
formula (1.8) prevents it from being so. We can get rid of this factor 
by taking the convolution square root. Let us recall that for a 
distribution $\mu$ its convolution square root $\mu^{\bplus 1/2}$ 
(in algebraic sense) is uniquely determined by 
\[
\Bigl(\mu^{\bplus 1/2}\Bigr)\,\bplus\, \Bigl(\mu^{\bplus 1/2}
\Bigr)=\mu
\qquad \mbox{or equivalently} \qquad
R(\mu^{\bplus 1/2})=\frac 12\cdot R(\mu).
\]
Thus we get the following striking corollary of Theorem 1.2.

$\ $

{\bf 1.11 Corollary.} The operation
\[
(\nu_1,\nu_2)\mapsto \lb \nu_1\Box\nu_2\rb^{\bplus 1/2}
\]
on distributions is associative.

$\ $

{\bf 1.12 Remarks.} 

1) This corollary raises the question, whether 
the operation $\lb\nu_1\bbox\nu_2\rb^{ \bplus 1/2}$ also makes sense 
in a $C^*$-context, i.e. whether the convolution square root of the 
free commutator of two probability measures is again a probability 
measure. The following counterexample shows that this is not the case
in general. 

Namely, let $p$ and $q$ be two free projections with
$\ff(q)=1/2$ and $0 < \ff(p)\leq 1/2$. We saw in Example 1.5.3.a 
that the distribution of $i(pq-qp)$ is the probability measure
\[
\mu_\xi \ = \ \sqrt{1- \xi} \ \delta_0 \ + \ \frac 1{\pi |t|}
\sqrt{\frac{4t^2-1+ \xi}{1-4t^2}} \ 
\chi_{\lb -1/2,-\sqrt{1- \xi}/2\rb\cup\lb\sqrt{1- \xi}/2,1/2\rb}dt,
\]
where $\xi\egdef 4\ff(p)(1-\ff(p))\in ( 0,1\rb$.
We claim now that $\mu_\xi^{\bplus 1/2}$ is not positive for
$\xi$ sufficiently close to 1. Let us assume, by contradiction, 
that there exist probability measures $\nu_\xi$ with
\[
\mu_\xi=\nu_\xi\,\bplus\, \nu_\xi\qquad (0 < \xi\leq 1).
\]
For $\xi=1$ (i.e. $\ff(p)=1/2$) the Application 1.3 gives us
\[
\nu_1=\frac 12 (\delta_{-1/4}+\delta_{1/4}).
\]
For $0 < \xi< 1$, $\mu_\xi$ has an atom at 0.
By Theorem 7.4 of \cite{BV2}, this implies that there exist atoms
$\alpha$ and $\beta$ for $\nu_\xi$ such that $\alpha+\beta=0$ and
$\nu_\xi(\{\alpha\})+\nu_\xi(\{\beta\})>1$.
Since $\mu_\xi$, and hence also $\nu_\xi$, is even we
can conclude that $\alpha=\beta=0$, i.e. $\nu_\xi$ has, for
all $0 < \xi<1$, an atom at 0 with $\nu_\xi(\{0\})>1/2$.
Consider now the mapping
$\xi\mapsto\nu_\xi$ for $\xi \in ( 0,1\rb$.
Clearly, all moments of $\nu_{\xi}$ are continuous in $\xi$.
Since we can find a compact interval which supports all $\nu_\xi$
(by the fact that the same is true for all $\mu_\xi$
and by invoking Lemma 3.1 of \cite{V2}), this implies that
$\xi\mapsto\nu_\xi$ is weakly continuous.
But this is in contradiction with the fact that for $\xi<1$ we
have mass at least $1/2$ at 0 whereas for $\xi=1$ all mass sits
at $\pm 1/4$.
Thus we can conclude that at least for $\xi$ sufficiently close
to 1 there exist no probability measure $\nu_\xi$ with
$\mu_\xi=\nu_\xi\,\bplus\,\nu_\xi$.

2) Theorem 1.2 shows that the main information needed for
calculating the free commutator consists of just the even free
cumulants. Thus it is natural to define a mapping 
$\mu\to Even(\mu)$ from all distributions to even distributions
by the requirement that $Even(\mu)$ is even (i.e. its 
odd free cumulants are zero) and that $Even(\mu)$ 
has the same even free cumulants as $\mu$:
\[
\lb R(Even(\mu))\rb(z)=\lb \ER(\mu)\rb (z^2).
\]
This raises the question whether this mapping also makes sense
on the $C^*$-level, i.e. whether it preserves positivity.
Using the same kind of arguing as before one can convince
oneself that this is not true in general.

Namely, according to Application 1.3, $\mu_{\xi}$ appearing above
is essentially (up to rescaling) $\mu_p\,\bplus\, \mu_{-p}$. But
this implies that $\nu_{\xi}$ is up to rescaling equal to
$Even(\mu_p)$. Since $\nu_{\xi}$ is not positive for $\xi$
sufficiently close to 1, positivity fails also for $Even(\mu_p)$ with
$\ff(p)$ sufficiently near to $1/2$.

$\ $

Whereas $\lb \nu_1\Box\nu_2\rb^{\bplus 1/2}$ does not make sense
for probability measures in general, it is again the case of even
probability measures which has a special status in this context.
As we have seen in Corollary 1.6, our results
take for even random variables an especially 
simple form if we formulate them in terms of the square
of the variables. In this spirit we will now make a transition from 
even probability measures to probability measures which have their
support on the positive real line, by just taking the square of the
corresponding random variables. We will denote by
$Q(\mu)$ the image of $\mu$ under the transformation
$Q:t\mapsto t^2$, i.e. we have $Q(\mu_a)=\mu_{a^2}$.

$\ $

{\bf 1.13 Corollary.} 

1) Let $\nu_1$ and $\nu_2$ be even probability 
measures. Then $\{ \nu_1\Box\nu_2\}^{\bplus 1/2}=
\lb \nu_1\Box\nu_2\rb^{\bplus 1/2}$ is an even probability measure, too.

2) The map $Q:t\mapsto t^2$ transports the operation
$(\nu_1,\nu_2)\mapsto \{ \nu_1\Box\nu_2\}^{\bplus 1/2}=
\lb\nu_1\bbox\nu_2\rb^{\bplus 1/2}$ on even probability measures
to the multiplicative free convolution
$(\nu_1,\nu_2)\mapsto \nu_1\,\btimes\, \nu_2$ on positively supported 
probability measures;
i.e., we have:
\begin{equation}
Q(\{ \nu_1\Box\nu_2\}^{\bplus 1/2})=
Q(\lb\nu_1\bbox\nu_2\rb^{\bplus 1/2})
=Q(\nu_1)\,\btimes\, Q(\nu_2).
\end{equation}
\newline

$\ $

Note that, although we will prove Corollary 1.13 only for 
probability measures with compact support, the continuity
properties of freeness (in particular Theorem 4.11 of
\cite{BV}) ensure that it is also true for arbitrary
even probability measures on $\R$. We next present an example
which uses this unbounded version of Corollary 1.13.

$\ $

{\bf 1.14 Example.} Let $a$ and $b$ be semicircular variables 
of radius 2 which are free. Then the free commutator of $a$ and $1/b$ 
has a Cauchy distribution of parameter $\gamma=2$, i.e.
\begin{equation}
\mu_{i(ab^{-1}-b^{-1}a)}=\frac \gamma{\pi(t^2+\gamma^2)}dt.
\end{equation}
This comes out from the Corollary 1.13 combined with the
results of \cite{Bia2} on free stable distributions, as follows.
The distribution of $1/b^2$ is the free stable distribution $\nu_\alpha$ 
of index $\alpha=1/2$, and on the other hand the distribution of
$a^2$ can be written as $\check \nu_{1/2}$ where $\check \nu$ denotes 
the image of $\nu$ by the map $t\mapsto 1/t$. The multiplicative free 
convolution $\nu_\alpha\,\btimes\,\check\nu_\alpha$ is given explicitly 
in Proposition 8 of \cite{Bia2} (comp. also \cite{Bia1}) in the form
\[
\nu_{1/2}\,\btimes\,\check\nu_{1/2}=\frac 1\pi \frac{t^{-1/2}}
{t+1}\chi_{(0,\infty)}dt,
\]
which turns out to be $Q^{-1}$ of the Cauchy distribution of parameter 
$\tilde\gamma=1$. Since the additive free convolution of a Cauchy 
distribution of parameter $\tilde\gamma$ with itself is a Cauchy distribution
of parameter $\gamma=2\tilde\gamma$, we obtain, by Corollary 1.13, 
that the distribution of $i(ab^{-1}-b^{-1}a)$ is Cauchy
of parameter $\gamma=2$.

$\ $

Let us now come back to the general structure of the free
commutator $\lb \nu_1\bbox\nu_2\rb$. Since $\freestar$
is commutative, Theorem 1.2 implies in particular that
also the free commutator is
commutative in its two variables. But the special form of our
formula (1.8) in Theorem 1.2 implies that we even have a form
of higher order commutativity.

$\ $

{\bf 1.15 Notations.}

1) A {\it commutator expression $f(\nu_1,\dots,\nu_n)$ 
of $n$ arguments} (for some 
$n\geq 1$) is an expression in which we take in some iterated 
way $n$ commutators of the distributions $\nu_1,\dots,\nu_n$.
More formally, we define the class of commutator expressions in
the following recursive way:
\newline
i) $f(\nu_1)=\nu_1$ is a commutator expression of one argument;
\newline
ii) $f(\nu_1,\dots,\nu_n)=\lb f_1(\nu_1,\dots,\nu_k)\bbox
f_2(\nu_{k+1},\dots,\nu_n)\rb$ is a commutator expression 
of $n$ arguments if, with
some $1\leq k< n$, $f_1$ and $f_2$ are commutator expressions
of $k$ and $n-k$ arguments, respectively.

2) To each commutator expression $f$ 
of $n$ arguments we assign a {\it depth vector} $(d_1,\dots,d_n)$
in the following recursive way:
\newline
i) $f(\nu_1)=\nu_1$ has depth $(0)$.
\newline
ii) If 
$f(\nu_1,\dots,\nu_n)=\lb f_1(\nu_1,\dots,\nu_k)\bbox
f_2(\nu_{k+1},\dots,\nu_n)\rb$ and if $f_1$ has depth 
$(d_1,\dots,d_k)$ and $f_2$ has depth $(d_{k+1},\dots,d_n)$,
then $f$ has depth $(d_1+1,\dots,d_k+1,d_{k+1}+1,\dots,d_n+1)$

$\ $

Intuitively, the depth vector of a commutator expression counts how many 
brackets we have to cross in order to reach the arguments inside the 
expression. For example, the expression
$\lb\, \lb\nu_1\bbox\nu_2\rb\bbox\lb\nu_3\bbox\nu_4\rb\,\rb$
has depth $(2,2,2,2)$, whereas 
$\lb\,\lb\,\lb \nu_1\bbox\nu_2\rb\bbox\nu_3\rb\bbox\nu_4\rb$
has depth $(3,3,2,1)$.

The announced higher order commutativity of the free commutator
can now be stated in the form that iterated commutators are
commutative in such variables which have the same depth,
i.e. which have equal entries in the depth vector.
More generally, we can compare two different commutator expressions
if their depth vectors differ only by a permutation.

$\ $

{\bf 1.16 Corollary.}
Let $f$ and $\widehat f$ be two commutator expressions of
$n$ arguments with depth $(d_1,\dots,d_n)$ and
$(\widehat d_1,\dots,\widehat d_n)$, respectively.
If there exists a permutation $\tau\in S_n$ such that
$(d_1,\dots,d_n)=(\widehat d_{\tau(1)},\dots,\widehat d_{\tau(n)})$
then
\begin{equation}
f(\nu_{\tau(1)},\dots,\nu_{\tau(n)}) \ = \ \widehat f(\nu_1,\dots,\nu_n)
\qquad\mbox{for arbitrary distributions $\nu_1,\dots,\nu_n$;}
\end{equation}
in particular
\begin{equation}
f(\nu,\dots,\nu)=\widehat f(\nu,\dots,\nu)
\qquad\mbox{for all distributions $\nu$.}
\end{equation}

$\ $

A concrete example of this corollary is, with
$f=\widehat f=\lb\lb \,\cdot\,\bbox\,\cdot\,\rb\bbox\lb\,\cdot\,
\bbox\,\cdot\,\rb\rb$, that
\[
\lb\, \lb\nu_1\bbox\nu_2\rb\bbox\lb\nu_3\bbox\nu_4\rb\,\rb=   
\lb\, \lb\nu_1\bbox\nu_3\rb\bbox\lb\nu_2\bbox\nu_4\rb\,\rb   
\]
for any distributions $\nu_1,\dots,\nu_4$.
 
$\ $

{\bf 1.17 Application (continuation of 1.7).}
Let us denote the canonical
iterated commutator of $m$ arguments by $f_m$, i.e.
\[
f_1(\nu_1):=\nu_1,\qquad f_m(
\nu_1,\dots,\nu_m):=\lb f_{m-1}(\nu_1,\dots,\nu_{m-1})
\bbox \nu_m\rb \quad (m\geq 2).
\]
Moreover, for every distribution $\nu$, let us denote 
\[
\nu_m:=f_m(\nu,\dots,\nu),\qquad (m\geq 1).
\]
Then we have
\begin{equation}
(\nu_n)_m=(\nu_m)_n\qquad \mbox{for all $m,n\geq 1$.}
\end{equation}
Indeed, given $m,n\geq 1$, let us denote by
$(f_m)_n$ the commutator expression of $m\cdot n$ arguments
defined by
\[
(f_m)_n(\nu_{11},\dots,\nu_{1m},\dots,\nu_{n1},\dots,\nu_{nm})=
f_n(f_m(\nu_{11},\dots,\nu_{1m}),\dots,
f_m(\nu_{n1},\dots,\nu_{nm})).
\]
An easy depth-counting argument shows that 
$(f_m)_n$ and $(f_n)_m$ satisfy the hypothesis of Corollary 1.16;
but then the Eqn. (1.32) applied to this situation gives
exactly (1.33).

Assume now that, in addition, the variance of $\nu$ is equal to
$1/2$. Then we know, by Application 1.7, that the moment-wise
limit $\nu_\infty:=\lim_{m\to\infty}\nu_m$ exists. By 
fixing $n$ and letting $m\to\infty$ in (1.33) it follows that:
\begin{equation}
(\nu_n)_\infty=(\nu_\infty)_n\qquad\mbox{for all $n\geq 1$.}
\end{equation}
If we take for $\nu$ the symmetric Bernoulli distribution
$1/2(\delta_{-1/2}+\delta_{1/2})$, then $\nu_2$ is the arcsine 
distribution on $\lb -1/2,1/2\rb$ (see, e.g., Example 1.5.3.a), 
hence we get the connection between the examples a) and b)
in Application 1.7 via the Eqn. (1.34), considered for $n=2$.

$\ $

An important ingredient in our proof of Theorem 1.2 is the notion of an 
$R$-diagonal pair, which was introduced and studied in our previous
paper \cite{NS3}. The connection of this concept with our present 
problem comes from the fact that if $a$ and $b$ are free
and both even, then $(ab,ba)$ is an $R$-diagonal pair --
which means that the joint distribution $\mu_{ab,ba}$ has
a nice, computable (two-dimensional) $R$-transform. This gives a
definite and treatable relation between the elements $ab$
and $ba$ (even though, of course, they are not free). 
The $R$-transforms of $\mu_{i(ab-ba)}$ and $\mu_{ab+ba}$
can be recovered from the one of $\mu_{ab,ba}$ by formulas
which have as immediate consequence that actually $\mu_{i(ab-ba)}
=\mu_{ab+ba}$ (as stated in Proposition 1.10). We take the occasion to 
note that another consequence of the $R$-diagonality of $(ab,ba)$ is the 
following polar decomposition result.

$\ $

{\bf 1.18 Proposition.}
Let $\ncps$ be a non-commutative probability space and assume that
$\A$ is a von Neumann algebra and $\ff$ is a faithful normal state.
Consider $a,b\in\A$ such that $a$ and $b$ are both selfadjoint and even,
$a$ is free from $b$, and $\mbox{Ker }a= \mbox{Ker }b=\{0\}$. 
Then $ab$ has a polar decomposition $ab=up$ where $u$ is a Haar unitary 
and where $\{u,u^*\}$ is free from $p$.

$\ $

The concatenation of arguments proving 1.10 and 1.18 is 
presented at the end of Section 4 of the paper.

$\ $

Finally, let us describe how the paper is organized. In the next section, 
we will collect the preliminaries on the combinatorial description of
freeness which will be needed in the forthcoming sections. In Section 3, 
we will reduce the general case to the case where $a$ and $b$ have even
distribution. In Section 4, we will recall the definition and the relevant 
results from \cite{NS3} on $R$-diagonal pairs and show that
if $a$ and $b$ are free and both even, then $(ab,ba)$ is an $R$-diagonal 
pair. This result will be used in Section 5 to conclude the proof 
of Theorem 1.2 and of its Corollaries 1.4, 1.6, 1.13, 1.16.

$\ $

$\ $

$\ $

\setcounter{section}{2}
\setcounter{equation}{0}
{\large\bf 2. Preliminaries}

$\ $

We start by collecting together a few basic free probabilistic
terms (some of them already implicitly reviewed and used in the 
Introduction).

$\ $

{\bf 2.1 Definitions.} 

1) We will call {\em non-commutative probability space} a pair $\ncps$,
where $\A$ is a unital algebra over $\C$, and 
$\varphi : \A \rightarrow \C$ is a linear functional 
normalized by $\varphi (1) =1.$ If we require in addition that 
$\A$ is a $C^{*}$-algebra and $\varphi$ is positive, then
$\ncps$ is called a {\em $C^{*}$-probability space.} 
If $\ff$ is a trace (i.e. $\ff(ab)=\ff(ba)$ for all $a,b\in\A$),
then we call $\ncps$ a {\em tracial} non-commutative
($C^*$-) probability space.

2) The {\em joint distribution} of a family of elements 
$\aunun \in \A$, in the non-commutative probability space 
$\ncps$, is the linear functional 
$\mu_{\aunun} : \ncpol \rightarrow \C$ given by:
\begin{equation}
\left\{  \begin{array}{l}
\mu_{\aunun} (1) \ = 1,  \\
\mu_{\aunun} ( \Xiqxin ) \ = \ \varphi ( \aiqaik ) \ \ \mbox{ for }
n \geq 1 \mbox{ and } 1 \leq \iunun \leq m, 
\end{array}  \right.
\end{equation}
where $\ncpol$ is the algebra of polynomials in $m$ non-commuting 
indeterminates $\Xunun$. 

3) If we only have one element $a \in \A$, then  
its distribution $\mu_{a} : \C [X] \rightarrow \C$, and the element  
$a$ itself are called {\em even} if $\mu_{a} (X^{n})$ (or equivalently,
$\varphi (a^{n})$) vanishes for all odd $n$.

4) A family of unital subalgebras $\A_{1}, \ldots , \A_{n} \subseteq \A$ 
is said to be {\em free} in the non-commutative probability space 
$\ncps$ if for every reduced word $w= a_{1} \cdots a_{k}$ made with
elements from the $\A_{i}$'s (i.e., $a_{1} \in \A_{i_{1}}, \ldots ,
a_{k} \in \A_{i_{k}}$ with $i_{1} \neq i_{2}, \ldots ,i_{k-1} \neq 
i_{k}$) we have the implication $\varphi (a_{1}) = \cdots =
\varphi (a_{k}) =0 \Rightarrow \varphi (w)=0.$ This definition
extends to subsets of $\A$, by putting $\X_{1} , \ldots , X_{n}
\subseteq \A$ to be free in $\ncps$ if and only if the unital
subalgebras generated by them are so.

5) Let $\ncps$ be a $C^{*}$-probability space. Then:

-- an element $a \in \A$ is called (standard) {\em semicircular} in 
$\ncps$ if $a=a^{*}$ and if $\mu_{a}$ is
the semicircle law $(2 \pi)^{-1} \sqrt{4-t^{2}} dt$ on $[-2,2];$

-- an element $a \in \A$ is called (standard) {\em circular} in
$\ncps$ if it is 
of the form $(a+ib)/ \sqrt{2},$ with $a,b$ free and semicircular;

-- an element $u \in \A$ is called a {\em Haar unitary} in $\ncps$
if it is a unitary and if $\varphi(u^n)=0$ for all 
$n\in\Z\backslash\{0\}$.

$\ $

We next review some combinatorial facts about non-crossing partitions.

\newpage

{\bf 2.2 Definitions.} 

1) If $\pi = \{ B_{1} , \ldots , B_{r} \}$
is a partition of $\nn$, 
then the equivalence relation on $\nn$ with equivalence classes
$B_{1} , \ldots , B_{r}$ will be denoted by $\ecpi$; the sets  
$B_{1} , \ldots , B_{r}$ are called the {\em blocks} of $\pi$. The 
number of elements in the set $B_k$ will be denoted by $|B_{k}|$.

2) A partition $\pi$ of $\nn$ is called {\em non-crossing} (\cite{K}) 
if for every $1 \leq i < j < i' < j' \leq n$ such that $i \ecpi i'$
and $j \ecpi j'$, it necessarily follows that
$i \ecpi j \ecpi i' \ecpi j'$. The set of all non-crossing 
partitions of $\nn$ will be denoted by $NC(n).$ 
By $NCE(n)$ we will denote the set of non-crossing partitions
$\pi\in NC(n)$ such that every block of $\pi$ has an even
number of elements. The complement $NC(n) \setminus NCE(n)$ will be
denoted by $NCO(n).$

3) On $NC(n)$ we will consider the partial order relation, called
{\em refinement order,} which is defined by
$\pi \leq \rho \ \ecdef$ each block of $\rho$ is a union of blocks 
of $\pi$. 

$\ $

{\bf 2.3 Definition (complementation maps on NC(n)).}
Let us now fix $n \geq 1$ and an $n$-tuple of 1's and 2's,
$\ee = ( \lunun ) \in \{ 1,2 \}^{n}.$ We describe here, following
\cite{NS3}, the construction of a map 
$C_{\ee} : NC(n) \rightarrow NC(n)$ which we will call 
{\em $\ee$-complementation map}, and which plays an important role
in the sequel.

Let us agree to call {\em $\ee$-insertion bead} a $2n$-tuple 
$( \Punun ; \Qunun )$ of points on a circle, constructed according to 
the following rules:

-- first, $\Punun$ are drawn around the circle, equidistant and in
clockwise order;

-- for every $1 \leq i \leq n$ such that $l_{i} =1,$ $Q_{i}$ is placed 
on the arc of circle going from $P_{i}$ to $P_{i+1}$ (clockwisely),
such that the length of the arc $P_{i}Q_{i}$ is 1/3 of the length of
the arc $P_{i}P_{i+1};$

-- for every $1 \leq i \leq n$ such that $l_{i} =2,$ $Q_{i}$ is placed 
on the arc of circle going from $P_{i}$ to $P_{i-1}$ (counter-clockwisely),
such that the length of the arc $P_{i}Q_{i}$ is 1/3 of the length of
the arc $P_{i}P_{i-1}.$
\newline
In the description of these rules, $P_{i+1}$ and $P_{i-1}$ are 
considered modulo $n$ (i.e. $P_{n+1} := P_{1}$ and $P_{0} := P_{n}$). 

Having fixed an $\ee$-insertion bead $( \Punun ; \Qunun )$, the
$\ee$-complement $\rho = C_{\ee} ( \pi ) \in NC(n)$ of a given partition
$\pi \in NC(n)$ is defined via the rule that for every $1 \leq i,j \leq n$
we have: $i \ecrho j \ \ecdef$ the line segment $Q_{i}Q_{j}$ does not
intersect any of the line segments $P_{h}P_{k}$, with $h \ecpi k.$

$\ $

To give a concrete example: for $n=5,$ $\ee = (1,1,2,2,1)$ and
$\pi = \{ \{ 1,2 \} , \{ 3,4,5 \} \}$ we have $C_{\ee} ( \pi )$ =
$\{ \{ 1 \} , \{ 2,3,5 \} , \{ 4 \} \}$, as shown by Figure 1.

\begin{figure}[t]

\setlength{\unitlength}{0.00066700in}%
\begingroup\makeatletter\ifx\SetFigFont\undefined
% extract first six characters in \fmtname
\def\x#1#2#3#4#5#6#7\relax{\def\x{#1#2#3#4#5#6}}%
\expandafter\x\fmtname xxxxxx\relax \def\y{splain}%
\ifx\x\y   % LaTeX or SliTeX?
\gdef\SetFigFont#1#2#3{%
  \ifnum #1<17\tiny\else \ifnum #1<20\small\else
  \ifnum #1<24\normalsize\else \ifnum #1<29\large\else
  \ifnum #1<34\Large\else \ifnum #1<41\LARGE\else
     \huge\fi\fi\fi\fi\fi\fi
  \csname #3\endcsname}%
\else
\gdef\SetFigFont#1#2#3{\begingroup
  \count@#1\relax \ifnum 25<\count@\count@25\fi
  \def\x{\endgroup\@setsize\SetFigFont{#2pt}}%
  \expandafter\x
    \csname \romannumeral\the\count@ pt\expandafter\endcsname
    \csname @\romannumeral\the\count@ pt\endcsname
  \csname #3\endcsname}%
\fi
\fi\endgroup
\begin{picture}(5700,5811)(3226,-6916)
\thicklines
\put(6001,-1636){\circle{300}}
\put(7726,-5686){\circle{300}}
\put(4351,-5686){\circle{300}}
\put(8251,-3211){\circle{300}}
\put(3751,-3211){\circle{300}}
\put(8401,-4111){\circle{300}}
\put(8176,-5011){\circle{300}}
\put(5251,-6286){\circle{300}}
\put(4201,-2386){\circle{300}}
\put(7126,-1786){\circle{300}}
\put(6151,-1711){\line( 4,-3){2100}}
\put(3901,-3286){\line( 5,-3){3750}}
\put(3826,-3361){\line( 1,-5){450}}
\put(4426,-5686){\line( 1, 0){3225}}
\put(8326,-3361){\line( 1,-6){ 75}}
\put(7651,-5761){\line(-5,-3){750}}
\put(6901,-6211){\line(-5,-1){750}}
\put(6151,-6361){\line(-1, 0){750}}
\put(4201,-5611){\line(-3, 4){225}}
\put(3976,-5311){\line(-1, 2){150}}
\put(3826,-5011){\line(-2, 5){150}}
\put(3676,-4636){\line(-1, 4){ 75}}
\put(3601,-3811){\line( 0, 1){225}}
\multiput(3601,-3586)(3.40909,10.22727){23}{\makebox(8.3333,12.5000){\SetFigFont{7}{8.4}{rm}.}}
\put(3751,-3061){\line( 1, 2){150}}
\put(3901,-2761){\line( 3, 4){225}}
\put(5926,-1561){\makebox(8.3333,12.5000){\SetFigFont{10}{12}{rm}.}}
\multiput(4951,-1786)(10.22727,3.40909){23}{\makebox(8.3333,12.5000){\SetFigFont{7}{8.4}{rm}.}}
\multiput(5176,-1711)(10.22727,3.40909){23}{\makebox(8.3333,12.5000){\SetFigFont{7}{8.4}{rm}.}}
\put(8401,-3961){\line( 0, 1){150}}
\put(5401,-1636){\line( 6, 1){450}}
\put(8401,-4261){\line(-1,-6){ 75}}
\multiput(8326,-4711)(-3.40909,-10.22727){23}{\makebox(8.3333,12.5000){\SetFigFont{7}{8.4}{rm}.}}
\put(3601,-3811){\makebox(8.3333,12.5000){\SetFigFont{10}{12}{rm}.}}
\put(3601,-4411){\line( 0, 1){225}}
\put(3601,-4186){\line( 0, 1){375}}
\put(4276,-2236){\line( 1, 1){225}}
\put(4501,-2011){\line( 2, 1){300}}
\multiput(4801,-1861)(9.37500,4.68750){17}{\makebox(8.3333,12.5000){\SetFigFont{7}{8.4}{rm}.}}
\put(6151,-1561){\line( 1, 0){375}}
\put(6526,-1561){\line( 3,-1){450}}
\put(7276,-1861){\line( 4,-3){300}}
\put(7576,-2086){\line( 1,-1){225}}
\put(7801,-2311){\line( 3,-4){225}}
\put(8026,-2611){\line( 1,-2){150}}
\multiput(8176,-2911)(3.40909,-10.22727){23}{\makebox(8.3333,12.5000){\SetFigFont{7}{8.4}{rm}.}}
\put(5101,-6286){\line(-2, 1){300}}
\multiput(4801,-6136)(-9.00000,6.00000){26}{\makebox(8.3333,12.5000){\SetFigFont{7}{8.4}{rm}.}}
\multiput(4576,-5986)(-7.50000,7.50000){21}{\makebox(8.3333,12.5000){\SetFigFont{7}{8.4}{rm}.}}
\put(7876,-5611){\line( 2, 3){300}}
\multiput(4201,-2536)(449.95663,-299.97109){9}{\line( 3,-2){225.347}}
\multiput(8101,-4861)(96.88521,484.42607){2}{\line( 1, 5){ 53.115}}
\multiput(8326,-4036)(-484.19220,193.67688){9}{\line(-5, 2){251.462}}
\put(5851,-1261){\makebox(0,0)[lb]{\smash{\SetFigFont{14}{16.8}{bf}P}}}
\put(8326,-2986){\makebox(0,0)[lb]{\smash{\SetFigFont{14}{16.8}{bf}P}}}
\put(3976,-6136){\makebox(0,0)[lb]{\smash{\SetFigFont{14}{16.8}{bf}P}}}
\put(3226,-3136){\makebox(0,0)[lb]{\smash{\SetFigFont{14}{16.8}{bf}P}}}
\put(8701,-4261){\makebox(0,0)[lb]{\smash{\SetFigFont{14}{16.8}{bf}Q}}}
\put(8551,-5311){\makebox(0,0)[lb]{\smash{\SetFigFont{14}{16.8}{bf}Q}}}
\put(7726,-6211){\makebox(0,0)[lb]{\smash{\SetFigFont{14}{16.8}{bf}P}}}
\put(4951,-6736){\makebox(0,0)[lb]{\smash{\SetFigFont{14}{16.8}{bf}Q}}}
\put(3601,-2161){\makebox(0,0)[lb]{\smash{\SetFigFont{14}{16.8}{bf}Q}}}
\put(5176,-6886){\makebox(0,0)[lb]{\smash{\SetFigFont{11}{13.2}{bf}4}}}
\put(4126,-6286){\makebox(0,0)[lb]{\smash{\SetFigFont{11}{13.2}{bf}4}}}
\put(3376,-3286){\makebox(0,0)[lb]{\smash{\SetFigFont{11}{13.2}{bf}5}}}
\put(3826,-2311){\makebox(0,0)[lb]{\smash{\SetFigFont{11}{13.2}{bf}5}}}
\put(6001,-1411){\makebox(0,0)[lb]{\smash{\SetFigFont{11}{13.2}{bf}1}}}
\put(8476,-3136){\makebox(0,0)[lb]{\smash{\SetFigFont{11}{13.2}{bf}2}}}
\put(8926,-4411){\makebox(0,0)[lb]{\smash{\SetFigFont{11}{13.2}{bf}2}}}
\put(8776,-5461){\makebox(0,0)[lb]{\smash{\SetFigFont{11}{13.2}{bf}3}}}
\put(7876,-6361){\makebox(0,0)[lb]{\smash{\SetFigFont{11}{13.2}{bf}3}}}
\put(6826,-1411){\makebox(0,0)[lb]{\smash{\SetFigFont{14}{16.8}{bf}Q}}}
\put(7051,-1561){\makebox(0,0)[lb]{\smash{\SetFigFont{11}{13.2}{bf}1}}}
\end{picture}

\end{figure}

\newpage

\[
\mbox{\bf{Figure 1.}}
\]

We shall be interested in the combinatorics of the $\ee$-complementation
map in two distinct situations. On one hand, we shall consider the case
when the sequence $\ee$ contains an equal number of 1's and 2's. Two
properties of $C_{\ee}$ for such $\ee$ are stated next.

$\ $

{\bf 2.4 Definition and Proposition.} Let $n \geq 1$ and 
$\ee = ( \lunun ) \in \{ 1,2 \}^{n}$ be such that $n$ is even and
$| \{ 1 \leq i \leq n \ | \ l_{i} =1 \} |$ =
$| \{ 1 \leq i \leq n \ | \ l_{i} =2 \} |$ = $n/2.$ We shall say that
a partition $\pi\in NC(n)$ is {\em $\ee$-alternating}
if for every block $B=\{i_1<i_2<\dots <i_k\}$ of $\pi$ we have
that $l_{i_1}\not=l_{i_2},\dots,l_{i_{k-1}}\not=l_{i_k},
l_{i_k}\not=l_{i_1}$. (Obviously, any such $\pi$ must be in $NCE(n),$
i.e. all its blocks must have even cardinality.) Given $\pi \in NC(n),$ 
the following hold:

1) (\cite{NS3}, Prop.7.7) If $\pi$ is $\ee$-alternating, then 
$C_\ee(\pi)$ is $\ee$-alternating, too.

2) (\cite{NS3}, Prop.8.11) If both $\pi$ and $C_{\ee} ( \pi )$
are in $NCE(n),$ then $\pi$ is $\ee$-alternating.

$\ $

At the other extreme, we shall be interested in the case when 
$\ee = (1,1, \ldots , 1).$

$\ $

{\bf 2.5 Definition.}  If $\ee = (1,1, \ldots ,1) \in \{ 1,2 \}^{n},$
then $C_{\ee}$ will be denoted by $K$, and is called the
{\em Kreweras complementation map} on $NC(n).$

$\ $

The map $K: NC(n) \rightarrow NC(n)$ was introduced in \cite{K}.
It turns out to be a bijection, and in fact an anti-automorphism
of $(NC(n), \leq )$ (i.e. $\pi \leq \rho \Leftrightarrow K( \pi )
\geq K( \rho )$ ). The inverse of $K$ is the $\ee$-complementation
corresponding to $\ee = (2,2, \ldots ,2).$

Recall that in the presentation of results made in Section 1, an
important role was played by the operation of combinatorial convolution
$\freestar$ (see e.g. Eqn. (1.8) of Theorem 1.2). By using the 
Kreweras complementation map, this operation can be defined as follows.

$\ $

{\bf 2.6 Definitions.}

1) Let $\Theta$ be the set of formal power series of the form 
$f(z) = \infsum \alpha_{n} z^{n},$ with $\alpha_{1}, \alpha_{2},
\alpha_{3}, \ldots  \in \C$. For $f \in \Theta$ and $n \geq 1,$ 
we shall use the notation $\coefn (f)$ for the coefficient of order $n$
of $f$. Moreover, for $f \in \Theta$, $n \geq 1$ and 
$\pi = \{ \Bunur \} \in NC(n),$ we shall use the notation
\begin{equation}
\coefnpi (f) \ := \ \prod_{i=1}^{r}  \ [ \mbox{coef } ( |B_{i}| ) ] (f).
\end{equation}

2) $\freestar$ is the binary operation on $\Theta$ which is defined by
the prescription:
\begin{equation}
\coefn ( f \  \freestar \  g) \ = \ \sum_{\pi \in NC(n)} \ 
\coefnpi (f) \cdot [ \mbox{coef } (n); K( \pi )] (g),
\end{equation}
for every $f,g \in \Theta$ and $n \geq 1.$

$\ $

For the interpretation of $\freestar$ which justifies the name of
`combinatorial convolution' (or, to be more precise: `convolution of
multiplicative functions on non-crossing partitions') we refer to
\cite{S1}, \cite{NS1}.

The operation $\freestar$ defined in 2.6 is associative, commutative, 
and has the series $Id(z)=z$ as a unit (see e.g. \cite{NS1}, Section 1.4).
An important role in the considerations related to $\freestar$ is played
by the series 
\begin{equation}
Zeta (z) \ := \ \infsum z^{n}, \ \ \ 
Moeb (z) \ := \ \infsum (-1)^{n+1} \frac{(2n-2)!}{(n-1)!n!} z^{n},
\end{equation}
called the Zeta and Moebius series, respectively.
$Zeta$ and $Moeb$ are inverse to each other with respect to $\freestar$.

Although chronologically the $R$-transform preceded the 
$\freestar$-operation, it is convenient to define it here in the 
following way.

$\ $

{\bf 2.7 Definition.} If $\mu$ is a distribution (in the algebraic sense
of Notations 1.1), then its {\em moment series} $M( \mu ) \in \Theta$ is
\begin{equation}
[ M( \mu ) ] (z) \ := \ \infsum \mu (X^{n}) z^{n},
\end{equation}
and its {\em $R$-transform} $R( \mu ) \in \Theta$ is
\begin{equation}
R( \mu ) \ := \ M( \mu )  \ \freestar \ Moeb.
\end{equation}
The coefficients of $R( \mu )$ are called the {\em free cumulants} 
of $\mu$.

$\ $

{\bf 2.8 Remark.} The Equation (2.6) can also be stated in the
equivalent form
\begin{equation}
M(\mu) \ = \ R(\mu)\,\freestar\, Zeta,
\end{equation}
which we will call the {\em moment-cumulant formula}. By identifying
the coefficient of order $n$ in (2.7), and by recalling how $\freestar$
is defined (Eqn. (2.3)), we get the explicit reformulation:
\begin{equation}
\mu(X^n)=\sum_{\pi=\{B_1,\dots,B_r\}\in NC(n)}
\alpha_{\vert B_1\vert}\dots \alpha_{\vert B_r\vert}\qquad (n\geq 1),
\end{equation}
where $( \alpha_{n} )_{n=1}^{\infty}$ are the free cumulants of $\mu$.

>From (2.8) (or (2.6-7)) it is clear that $\mu \rightarrow R( \mu )$ is a 
bijection from the set of distributions onto $\Theta$. Another fact
immediately implied by (2.8) is that a distribution $\mu$ is even (in 
the sense of Definition 2.1.3) if and only if all its odd free cumulants 
are equal to zero.

$\ $

{\bf 2.9 Remark.} The direct connection between the operation $\freestar$
and freeness is provided by the fact, already mentioned in the Introduction,
that:
\begin{equation}
R( \mu_{ab} ) \ = \ R( \mu_{a} ) \ \freestar \ R( \mu_{b} )
\end{equation}
whenever $a,b$ are free in some non-commutative probability space $\ncps$.
For the proof of (2.9), see \cite{NS1}, Section 3.5. This equation can be 
reformulated as follows: first, we $\freestar$-operate with $Zeta$ on both
its sides and thus (by also taking (2.7) into account) we turn it into:
\begin{equation}
M( \mu_{ab} ) \ = \ R( \mu_{a} ) \ \freestar \ M( \mu_{b} );
\end{equation}
then by taking the coefficient of order $n$ in (2.10), we obtain the formula:
\begin{equation}
\varphi ( (ab)^{n} ) \ = \ \sum_{\pi \in NC(n)} 
\coefnpi (R( \mu_{a} )) \cdot [ \mbox{ coef} (n);K( \pi )](M( \mu_{b} )),
\end{equation}
valid for $a,b$ free in $\ncps$, and $n \geq 1.$

In the Sections 3 and 4 below we will need the generalization of (2.11) to 
the situation when $K$ is replaced by an arbitrary $\ee$-complementation 
map $C_{\ee}$; this is given in the next proposition, and is a particular 
case of Proposition 7.3 in \cite{NS3}.

$\ $

{\bf 2.10 Proposition (\cite{NS3}).} Let $a,b$ be free in the non-commutative 
probability space $\ncps$, and let us denote $ab=:x_{1},$ $ba=:x_{2}.$ Then 
for every $n \geq 1$ and $\ee = ( \lunun ) \in \{ 1,2 \}^{n}$ we have 
\begin{equation}
\varphi ( x_{l_{1}} \cdots x_{l_{n}} ) \ = \ \sum_{\pi \in NC(n)} \ 
\coefnpi (R( \mu_{a} )) \cdot [ \mbox{coef } (n); C_{\ee} ( \pi )]
(M( \mu_{b} )).
\end{equation}

$\ $

Even though all the results in the Introduction are concerning 1-dimensional
distributions, the key fact in their proofs will be that we will often 
work in 2 variables. We thus conclude this section by reviewing the facts
that will be needed about the 2-dimensional $R$-transform; this goes in 
parallel with the path taken in the Definitions 2.6, 2.7.

$\ $

{\bf 2.11 Definitions.}

1) Let $\Theta_{2}$ be the set of formal power series without constant
coefficient, in 2 non-commuting variables $z_{1},z_{2}.$
An element of $\Theta_{2}$ is thus a series of the form
\begin{equation}
f( z_{1} , z_{2} ) \ = \ \sum_{n=1}^{\infty} \ \ 
\sum_{( \iunun ) \in \{ 1,2 \}^{n}} \ \alpha_{( \iunun )} \ziqzin ,
\end{equation} 
where the $\alpha_{( \iunun )}$'s are some complex coefficients.
For $f \in \Theta_{2}$, $n \geq 1$ and $( \iunun ) \in \{ 1,2 \}^{n}$ 
we shall use the notation $\coefnnn (f)$ for the coefficient of 
$z_{i_{1}} \cdots z_{i_{n}}$ in $f$. Moreover, for $f \in \Theta_{2}$, 
$n \geq 1,$ $( \iunun ) \in \{ 1,2 \}^{n}$ and 
$\pi = \{ \Bunur \} \in NC(n)$ we shall use the notation
\begin{equation}
\coefnnnpi (f) \ := \ \prod_{j=1}^{r} \ [ \mbox{coef } 
( \iunun \ | \ B_{j} ) ] (f),
\end{equation}
where in the right-hand side of (2.14) we write
\begin{equation}
( \iunun \ | \ B) \ := \ ( i_{h_{1}}, \ldots , i_{h_{m}} )
\end{equation}
whenever $B = \{ h_{1} < \cdots < h_{m} \}$ is a non-void subset of $\nn$.

2) We denote by $\freestar_{2}$ the binary operation on $\Theta_{2}$ which 
is defined by the prescription:
\begin{equation}
\coefnnn ( f \ \freestar_{2} \ g) \ = \ \sum_{\pi \in NC(n)} 
\coefnnnpi (f) \cdot [ \mbox{coef } ( \iunun ); K( \pi )] (g),
\end{equation}
for every $f,g \in \Theta_{2},$ $n \geq 1$ and $( \iunun ) \in \{ 1,2 \}^{n}.$
Then $\freestar_{2}$ is associative (but not commutative), and has 
the series $Sum(z_{1},z_{2})=z_{1}+z_{2}$ as a unit. The series 
\[
Zeta_{2} (z_{1},z_{2}) \ := \ \infsum \ \sum_{ ( \iunun ) \in 
\{ 1,2 \}^{n} } \ziqzin ,
\]
\begin{equation}
\end{equation}
\[
Moeb_{2} (z_{1},z_{2}) \ := \ \infsum \ \sum_{ ( \iunun ) \in
\{ 1,2 \}^{n} } \ (-1)^{n+1} \frac{(2n-2)!}{(n-1)!n!} \ziqzin ,
\]
are called the 2-variable Zeta and Moebius series; they are inverse
to each other with respect to $\freestar_{2},$ and they are central (i.e.
$f \ \freestar \ Zeta_{2} = Zeta_{2} \ \freestar \ f$ for every 
$f \in \Theta_{2},$ and similarly for $Moeb$). For the proofs of all these 
facts, see \cite{NS2}, Section 3.

3) If $\mu : \C \langle X_{1},X_{2} \rangle \rightarrow \C$ is a linear
functional normalized by $\mu (1)=1,$ then its {\em moment series} 
$M( \mu ) \in \Theta_{2}$ is
\begin{equation}
[ M( \mu ) ] (z_{1},z_{2}) \ := \ \infsum \
\sum_{( \iunun ) \in \{ 1,2 \}^{n}} \  \mu ( \Xiqxin ) \ziqzin ,
\end{equation}
and its {\em $R$-transform} $R( \mu ) \in \Theta_{2}$ is
\begin{equation}
R( \mu ) \ := \ M( \mu )  \ \freestar_{2} \ Moeb_{2}.
\end{equation}
The coefficients of $R( \mu )$ are called the the {\em free cumulants} 
of $\mu$.

$\ $

{\bf 2.12 Remark.} Same as in the 1-dimensional case, Eqn. (2.19) can be 
stated in the equivalent form
\begin{equation}
M( \mu ) \ = \ R( \mu ) \ \freestar_{2} \ Zeta_{2} ,
\end{equation}
called the {\em moment-cumulant formula} (in two variables). From 
(2.19--20) it follows that $\mu \rightarrow R( \mu )$ is a bijection from 
the set of normalized linear functionals on $\C \langle X_{1}, X_{2} \rangle$ 
onto $\Theta_{2}.$ Also, by identifying the coefficient of $\ziqzin$ on the 
two sides of (2.20), one can easily prove the following fact: let 
$\mu : \C \langle X_{1},X_{2} \rangle \rightarrow \C$ be a normalized linear 
functional, with free cumulants $( \alpha_{\iunun} )_{n; \iunun}$; then 
`$\mu ( \Xiqxin ) = 0$ for every odd $n$ and every $( \iunun ) 
\in \{ 1,2, \}^{n}$' is equivalent to `$\alpha_{\iunun} =0$ for every
odd $n$ and every $( \iunun ) \in \{ 1,2 \}^{n}$'.

We mention that the important Equation (2.9), relating $\freestar$ with
the multiplication of free random variables, extends to the multi-dimensional
case; in 2 variables we have: 
\begin{equation}
R( \mu_{a_{1}b_{1} , a_{2}b_{2}} ) \ = \ 
R( \mu_{a_{1} , a_{2} } )  \ \freestar_{2}  \ R( \mu_{b_{1} , b_{2} } ) ,
\end{equation}
whenever $\{ a_{1},a_{2} \}$ is free from $\{ b_{1}, b_{2} \}$ in a
non-commutative probability space $\ncps$. The multi-variable versions of
(2.9) have some interesting applications to freeness, see \cite{NS2}. 

We also state here, for future reference, a
`partly 2-variable' version of Proposition 2.10; this is still a particular 
case of Proposition 7.3 in \cite{NS3}.

$\ $

{\bf 2.13 Proposition (\cite{NS3}).} Let $a$ be free from $\{ b_{1},b_{2}
\}$ in the non-commutative probability space $\ncps$, and let us denote 
$ab_{1}=:y_{1},$ $b_{2}a=:y_{2}.$ Then for every
$n \geq 1$ and $\ee = ( \lunun ) \in \{ 1,2 \}^{n}$ we have 
\begin{equation}
\varphi ( y_{l_{1}} \cdots y_{l_{n}} ) \ = \
\sum_{\pi \in NC(n)}
\coefnpi (R( \mu_{a} )) \cdot [ \mbox{coef } ( \lunun );C_{\ee} ( \pi )]
(M( \mu_{b_{1},b_{2}} )).
\end{equation}

$\ $

The 2-dimensional $R$-transform comes into the free commutator
problem via the following fact, which is an immediate application of 
how the $R$-transform behaves under linear transformations (see \cite{N},
Section 5).

$\ $

{\bf 2.14 Proposition.} Let $\ncps$ be a non-commutaive probability space,
and let $a,b$ be in $\A$. Then we have the relations:
\begin{equation}
[ R( \mu_{i(ab-ba)} ) ] (z) \ = \ [ R( \mu_{ab,ba} ) ] (iz,-iz),
\end{equation}
and
\begin{equation}
[ R( \mu_{ab+ba} ) ] (z) \ = \ [ R( \mu_{ab,ba} ) ] (z,z).
\end{equation}

$\ $

$\ $

$\ $

\setcounter{section}{3}
\setcounter{equation}{0}
{\large\bf 3. Reduction to the case of even random variables} 

$\ $

The result which will enable us to make such a reduction 
is the following.

$\ $

{\bf 3.1 Theorem.} Let $\ncps$ and $\pstild$ be $\noncom$s, and let
$a,b \in \A$, $\atild , \btild \in \Atild$ be such that:

i) $a$ is free from $b$ in $\ncps$, $\atild$ is free from $\btild$ in
$\pstild$;

ii) $\RE(\mu_a)=\RE(\mu_{\atild})$ (with $\er$ as in Notations 1.1);

iii) $\RE(\mu_b)=\RE(\mu_{\btild})$.

Then the distribution of $ab-ba$ in $\ncps$ coincides with the one of
$\atild \btild - \btild \atild$ in $\pstild .$

$\ $

The way how Theorem 3.1 is used will be shown later on in the paper
(Section 5.4); for the moment we only concentrate on its proof. 
We start by making a further reduction in the statement of the theorem.

$\ $

{\bf 3.2 Lemma.} It suffices to prove Theorem 3.1 in the case when its
hypothesis (iii) is replaced by the stronger one that:
(iii$'$) $\mu_{b} = \mu_{\btild}$. 

$\ $

{\bf Proof.} Let us assume the Theorem 3.1 to be true when the hypothesis
(i)+(ii)+(iii$'$) are holding, and let us consider $\ncps$, $\pstild$,
$a,b \in \A$, $\atild , \btild \in \Atild$ which only satisfy 
(i)+(ii)+(iii). By using a free product construction, we can produce a
$\noncom$ $( \widehat{\A} , \widehat{\varphi} )$ and elements
$\widehat{a} , \widehat{b} \in \widehat{A}$ which are free with respect to
$\widehat{\varphi}$ and such that $\mu_{\widehat{a}} = \mu_{a},$
$\mu_{\widehat{b}} = \mu_{\widetilde{b}} .$ Then on one hand,
$\atild , \btild \in \Atild$ and
$\widehat{a} , \widehat{b} \in \widehat{A}$ satisfy (i)+(ii)+(iii$'$),
hence $\comtild$ has the same distribution as $\widehat{a} \widehat{b} -
\widehat{b} \widehat{a} ;$ on the other hand,
$b,a \in \A$ and
$\widehat{b} , \widehat{a} \in \widehat{A}$ also satisfy (i)+(ii)+(iii$'$),
hence $ba-ab$ has the same distribution as $\widehat{b} \widehat{a} -
\widehat{a} \widehat{b} .$
>From these two facts it is immediate that $\comm$ and $\comtild$ have 
identical distributions. {\bf QED}

$\ $

We shall prove the Theorem 3.1 (in the reduced form of 3.2) by exhibiting
an expression for the moments of $\comm$, where only the even free cumulants
of $a$ are involved. The first step towards doing this is the following.

$\ $

{\bf 3.3 Proposition.} Let $\ncps$ be a $\noncom$, and let $a,b \in \A$
be free with respect to $\varphi .$ Then for every $n \geq 1$ we have:
\begin{equation}
\varphi ( ( \comm )^{n} ) \ = \ 
\sum_{ \begin{array}{c}
{\scriptstyle \pi \in NC(n) }  \\
{\scriptstyle \ee \in \strings}
\end{array}  } \  (-1)^{d( \ee )} \
[ \mbox{coef } (n) ; \pi ] ( \rmua ) \cdot
[ \mbox{coef } (n) ; C_{\ee} ( \pi ) ] ( \mmub ),
\end{equation}
where for an $n$-tuple $\ee \in \strings$, $d ( \ee )$ denotes the number 
of components of $\ee$ that are equal to 2 (and the notations $C_{\ee},$ 
$\coefnpi$ are as established in Section 2).

$\ $

{\bf Proof.} Let us make the notations $ab =x_{1},$ $ba = x_{2},$ and for 
every $\ee \egdef ( \lunun ) \in \strings$ let us put 
$x_{\ee} \egdef x_{l_{1}} x_{l_{2}} \cdots x_{l_{n}} \in \A .$ It is then
obvious that 
\begin{equation}
\varphi ( ( \comm )^{n} ) \ = \ \varphi ( ( x_{1} - x_{2} )^{n} ) \ = \ 
\sum_{\ee \in \strings} \ (-1)^{d( \ee )} \varphi ( x_{\ee} ) .
\end{equation}
Equation (3.1) follows from (3.2) and the expression for $\varphi 
( x_{\ee} )$ as summation over non-crossing partitions which 
is provided by Eqn. (2.12) of Proposition 2.10. {\bf QED}

$\ $

The second step is to put in evidence some remarkable involutions on 
$NCO(n)$ and $NCO(n) \times \{ 1,2 \}^{n}.$

$\ $

{\bf 3.4 Notations.} 

1) Recall from Section 2.2 that we denote by $NCO(n)$ the set of 
non-crossing partitions $\pi$ of $\nn$ which have at least one block with
an odd number of elements. Given $\pi \in NCO(n),$ it is easy to show
that there exist blocks $B$ of $\pi$ having $|B|$ odd and 
min $B$, max $B$ of the same parity (see e.g. \cite{NS3}, Lemma 4.7);
among these blocks let us pick the one, call it $B_{o}$, which has the 
smallest value of min $B_{o}.$ It is convenient that (for the purposes of
this section only) we give a name to the interval [min $B_{o}$, max
$B_{o}$]; we shall call it the {\em twist interval} of $\pi$.

2) For $\pi \in NCO(n)$ we shall denote by $tw( \pi )$ (with $tw$ for
`twist') the partition of $\nn$ obtained in the following way. Let 
$[t_{0},t_{1}]$ be the twist interval of $\pi$, as defined above, and
let $\tau$ be the permutation of $\nn$ given by:
$\tau (i) = t_{0}+t_{1}-i,$ if $t_{0} \leq i \leq t_{1},$ and 
$\tau (i) = i$ otherwise. Then we put $tw( \pi )$ to be 
`$\tau ( \pi )$', i.e.:
\begin{equation}
tw( \pi ) \ := \ \{ \tau (B) \ | \ B \mbox{ block of } \pi \} .
\end{equation}
It is immediate that $tw( \pi )$ is also in $NCO(n)$, and has the same
twist interval $[t_{0},t_{1}]$ as $\pi$ itself. The latter fact has as 
consequence that $tw(tw( \pi )) = \pi$; i.e., $tw : NCO(n) \rightarrow
NCO(n)$ is an involution.

3) We need to adapt the twist $tw$ to the case when an $n$-tuple 
$\ee \in \{ 1,2 \}^{n}$ is considered at the same time with the partition
$\pi \in NCO(n).$ Thus, we shall denote by $\twtw$ the map from
$NCO(n) \times \{ 1,2 \}^{n}$ into itself which is defined by the formula:
\begin{equation}
\twtw ( \pi , ( \lunun )) \ = \ ( \ tw( \pi ), \ ( l_{1} ' , \ldots ,
l_{n} ' ) \ ), \ \ \ \pi \in NCO(n), \  \lunun \in \{ 1,2 \} ,
\end{equation}
where
\begin{equation}
l_{i} ' \ = \  \left\{  \begin{array}{lll}
3- l_{t_{0}+t_{1}-i} &  &  \mbox{ if } t_{0} \leq i \leq t_{1}  \\
l_{i}                &  &  \mbox{ otherwise, } 
\end{array}  \right.
\end{equation}
and where in (3.5) $t_{0}$ and $t_{1}$ are the end-points of the twist 
interval of $\pi$.

If $\twtw ( \pi , \ee ) =: ( tw ( \pi ), \ee ' )$ and
if $\twtw ( tw( \pi ) , \ee ' ) =: (  \pi , \ee '' )$, then from (3.5) and
the fact that $tw( \pi )$ has the same twist interval as $\pi$, we infer
that $\ee '' = \ee$. Hence $\twtw$ is an involution on
$NCO(n) \times \{ 1,2 \}^{n}.$

$\ $

{\bf 3.5 Lemma.} If $( \pi , \ee )$ is in $NCO(n) \times \{ 1,2 \}^{n}$
for some $n,$ and if $\twtw ( \pi , \ee ) =: ( \pi ' , \ee ' ),$ then
the partitions $C_{\ee} ( \pi )$ and $C_{\ee '} ( \pi ')$ have
the same block structure (i.e., for every $1 \leq m \leq n,$ the two
partitions have the same number of blocks with $m$ elements).

$\ $

{\bf Proof.} Let $[t_{0},t_{1}]$ be the twist interval of $\pi$, and
let $\tau$ be the permutation of $\nn$ defined by
$\tau (i) = t_{0}+t_{1}-i,$ if $t_{0} \leq i \leq t_{1},$ and 
$\tau (i) = i$ otherwise. We also denote $C_{\ee} ( \pi ) =: \rho$,
$C_{\ee '} ( \pi ' ) =: \rho '.$ We will show that $\rho = \tau 
( \rho ')$ (where $\tau ( \rho ')$ is defined as $\{ \tau (B) \ |$ 
$B$ block of $\rho ' \}$ ); this will of course imply that $\rho$ and
$\rho '$ have the same block structure.

It will actually suffice to show that 
\begin{equation}
\rho \geq \tau ( \rho ') \ \ \ \mbox{ (in the refinement order);}
\end{equation}
since at the moment it is not verified that $\tau ( \rho ')$ is 
non-crossing, we consider this inequality in the larger poset of all
the partitions of $\nn$.

Indeed, let us assume that (3.6) is proved. Then by applying this
inequality to $( \pi ', \ee ')$ instead of $( \pi , \ee )$, and by taking
into account that $\twtw ( \pi ' , \ee ') = ( \pi , \ee )$ and that
$\pi '$ has the same twist interval as $\pi$, we get: $\rho ' \geq
\tau ( \rho ).$ But $\tau^{2} = id,$ hence the application of $\tau$ to
the both sides of the latter inequality leads to $\tau ( \rho ') \geq \rho$,
and hence (3.6) must be an equality.

Due to the way how $\rho = C_{\ee} ( \pi )$ and $\rho ' =
C_{\ee '} ( \pi ')$ are concretely defined, the proof of (3.6) is an exercise
in elementary geometry. We consider an $\ee$-insertion bead 
$( \Punun ; Q_{1}, \ldots ,$
$Q_{n})$ (see Definition 2.3), and on the same circle we mark
the points $X_{0},X_{1},Y_{0},Y_{1}$ in the way shown in Figure 2; if the
length of the arc of circle between two consecutive $P_{i}$'s is denoted
by $\lambda$, then the length of the arcs $Y_{0}P_{t_{0}}$ and
$P_{t_{1}}Y_{1}$ is $\lambda /3$, and the one of the arcs $X_{0}Y_{0}$ and
$Y_{1}X_{1}$ is $\lambda / 6.$

\begin{figure}[t]

\setlength{\unitlength}{0.00066700in}%
\begingroup\makeatletter\ifx\SetFigFont\undefined
% extract first six characters in \fmtname
\def\x#1#2#3#4#5#6#7\relax{\def\x{#1#2#3#4#5#6}}%
\expandafter\x\fmtname xxxxxx\relax \def\y{splain}%
\ifx\x\y   % LaTeX or SliTeX?
\gdef\SetFigFont#1#2#3{%
  \ifnum #1<17\tiny\else \ifnum #1<20\small\else
  \ifnum #1<24\normalsize\else \ifnum #1<29\large\else
  \ifnum #1<34\Large\else \ifnum #1<41\LARGE\else
     \huge\fi\fi\fi\fi\fi\fi
  \csname #3\endcsname}%
\else
\gdef\SetFigFont#1#2#3{\begingroup
  \count@#1\relax \ifnum 25<\count@\count@25\fi
  \def\x{\endgroup\@setsize\SetFigFont{#2pt}}%
  \expandafter\x
    \csname \romannumeral\the\count@ pt\expandafter\endcsname
    \csname @\romannumeral\the\count@ pt\endcsname
  \csname #3\endcsname}%
\fi
\fi\endgroup
\begin{picture}(8497,6719)(2251,-7208)
\thicklines
\put(7126,-1786){\circle{300}}
\put(6226,-6361){\circle{300}}
\put(5401,-6361){\circle{336}}
\put(4801,-6136){\circle{336}}
\put(6301,-1561){\circle{336}}
\put(5701,-1561){\circle{336}}
\put(8326,-3361){\line( 1,-6){ 75}}
\put(4201,-5611){\line(-3, 4){225}}
\put(3976,-5311){\line(-1, 2){150}}
\put(3826,-5011){\line(-2, 5){150}}
\put(3676,-4636){\line(-1, 4){ 75}}
\put(3601,-3811){\line( 0, 1){225}}
\multiput(3601,-3586)(3.40909,10.22727){23}{\makebox(8.3333,12.5000){\SetFigFont{7}{8.4}{rm}.}}
\put(3751,-3061){\line( 1, 2){150}}
\put(3901,-2761){\line( 3, 4){225}}
\put(5926,-1561){\makebox(8.3333,12.5000){\SetFigFont{10}{12}{rm}.}}
\put(8401,-3961){\line( 0, 1){150}}
\put(8401,-4261){\line(-1,-6){ 75}}
\multiput(8326,-4711)(-3.40909,-10.22727){23}{\makebox(8.3333,12.5000){\SetFigFont{7}{8.4}{rm}.}}
\put(3601,-3811){\makebox(8.3333,12.5000){\SetFigFont{10}{12}{rm}.}}
\put(3601,-4411){\line( 0, 1){225}}
\put(3601,-4186){\line( 0, 1){375}}
\put(4276,-2236){\line( 1, 1){225}}
\put(4501,-2011){\line( 2, 1){300}}
\multiput(4801,-1861)(9.37500,4.68750){17}{\makebox(8.3333,12.5000){\SetFigFont{7}{8.4}{rm}.}}
\put(6526,-1561){\line( 3,-1){450}}
\put(7276,-1861){\line( 4,-3){300}}
\put(7576,-2086){\line( 1,-1){225}}
\put(7801,-2311){\line( 3,-4){225}}
\put(8026,-2611){\line( 1,-2){150}}
\multiput(8176,-2911)(3.40909,-10.22727){23}{\makebox(8.3333,12.5000){\SetFigFont{7}{8.4}{rm}.}}
\put(7876,-5611){\line( 2, 3){300}}
\put(7126,-1936){\line(-1,-5){900}}
\multiput(5851,-3961)(512.15845,-102.43169){10}{\line( 5,-1){265.574}}
\put(4426,-5836){\line(-1, 1){225}}
\put(3676,-3361){\line( 1, 5){ 75}}
\multiput(4126,-2461)(6.00000,9.00000){26}{\makebox(8.3333,12.5000){\SetFigFont{7}{8.4}{rm}.}}
\put(8401,-4261){\line( 0, 1){300}}
\multiput(8176,-5161)(3.40909,10.22727){23}{\makebox(8.3333,12.5000){\SetFigFont{7}{8.4}{rm}.}}
\multiput(8326,-3361)(-3.40909,10.22727){23}{\makebox(8.3333,12.5000){\SetFigFont{7}{8.4}{rm}.}}
\put(4426,-5836){\line( 4,-3){300}}
\put(4951,-6211){\line( 2,-1){300}}
\put(5551,-6436){\line( 1, 0){525}}
\put(6376,-6436){\line( 5, 1){375}}
\put(6751,-6361){\line( 4, 1){300}}
\multiput(7051,-6286)(10.22727,3.40909){23}{\makebox(8.3333,12.5000){\SetFigFont{7}{8.4}{rm}.}}
\multiput(7276,-6211)(9.00000,6.00000){26}{\makebox(8.3333,12.5000){\SetFigFont{7}{8.4}{rm}.}}
\put(7501,-6061){\line( 1, 1){225}}
\multiput(7726,-5836)(6.00000,9.00000){26}{\makebox(8.3333,12.5000){\SetFigFont{7}{8.4}{rm}.}}
\put(5851,-1561){\line( 1, 0){300}}
\put(4951,-1786){\line( 5, 2){375}}
\multiput(5326,-1636)(10.22727,3.40909){23}{\makebox(8.3333,12.5000){\SetFigFont{7}{8.4}{rm}.}}
\put(9751,-5161){\vector(-1,-4){  0}}
\put(9751,-5161){\vector( 1, 4){185.294}}
\multiput(5701,-1786)(-110.26932,-551.34658){8}{\line(-1,-5){ 53.115}}
\put(5326,-3886){\line( 6,-1){450}}
\multiput(4801,-6286)(-105.47527,-632.85165){2}{\line(-1,-6){ 44.525}}
\multiput(4651,-7186)(0.00000,8.82353){9}{\line( 0, 1){  4.412}}
\multiput(5776,-1411)(105.47527,632.85165){2}{\line( 1, 6){ 44.525}}
\put(5401,-1261){\makebox(0,0)[lb]{\smash{\SetFigFont{14}{16.8}{bf}X}}}
\put(6226,-1261){\makebox(0,0)[lb]{\smash{\SetFigFont{14}{16.8}{bf}Y}}}
\put(7051,-1411){\makebox(0,0)[lb]{\smash{\SetFigFont{14}{16.8}{bf}P}}}
\put(5251,-6886){\makebox(0,0)[lb]{\smash{\SetFigFont{14}{16.8}{bf}Y}}}
\put(5626,-1336){\makebox(0,0)[lb]{\smash{\SetFigFont{11}{13.2}{bf}0}}}
\put(6376,-1336){\makebox(0,0)[lb]{\smash{\SetFigFont{11}{13.2}{bf}0}}}
\put(5476,-6961){\makebox(0,0)[lb]{\smash{\SetFigFont{11}{13.2}{bf}1}}}
\put(6151,-6961){\makebox(0,0)[lb]{\smash{\SetFigFont{14}{16.8}{bf}P}}}
\put(6301,-7036){\makebox(0,0)[lb]{\smash{\SetFigFont{11}{13.2}{bf}t}}}
\put(6376,-7111){\makebox(0,0)[lb]{\smash{\SetFigFont{11}{13.2}{bf}1}}}
\put(7201,-1486){\makebox(0,0)[lb]{\smash{\SetFigFont{11}{13.2}{bf}t}}}
\put(7276,-1561){\makebox(0,0)[lb]{\smash{\SetFigFont{11}{13.2}{bf}0}}}
\put(9826,-4261){\makebox(0,0)[lb]{\smash{\SetFigFont{14}{16.8}{bf}A}}}
\put(9601,-5536){\makebox(0,0)[lb]{\smash{\SetFigFont{14}{16.8}{bf}T(A)}}}
\put(4351,-6511){\makebox(0,0)[lb]{\smash{\SetFigFont{14}{16.8}{bf}X}}}
\put(4576,-6586){\makebox(0,0)[lb]{\smash{\SetFigFont{11}{13.2}{bf}1}}}
\put(2251,-5161){\makebox(0,0)[lb]{\smash{\SetFigFont{14}{16.8}{bf}B=T(B)}}}
\end{picture}
 
\end{figure}

Let us record the remark that:
\begin{equation}
[ P_{h} , P_{k} ] \cap [ Y_{0} , Y_{1} ] \ = \emptyset \ =
[ P_{h} , P_{k} ] \cap [ X_{0} , X_{1} ] 
\end{equation}
for every $1 \leq h ,k \leq n$ such that $h \ecpi k$, and where
we use the notation $[A,B]$ for the line segment connecting the points
$A,B$ in the plane. The verification of (3.7) is easily made by using that
$t_{0}$ and $t_{1}$ are the min and max of the same block of $\pi$
(hence $\{t_{0} , \ldots , t_{1} \}$ and $\nn \setminus \{ t_{0}, \ldots ,
t_{1} \}$are both unions of blocks of $\pi$).
 
Let us next denote by ${\cal S}$ the open-half plane determined by
$X_{0}$ and $X_{1},$ which contains $Y_{0},Y_{1},P_{t_{0}},P_{t_{1}}.$
We shall denote by $T$ the bijective transformation of the plane 
which reflects the points of ${\cal S}$ in the mediator line of
$[X_{0},X_{1}]$, and leaves fixed the points in the complement of 
${\cal S}$ (see Figure 2). It is obvious that
\begin{equation}
T ( P_{i} ) \ = \ P_{\tau (i)} , \ \ \ 1 \leq i \leq n,
\end{equation}
where $\tau$ is the permutation mentioned in the first phrase of the proof.
We do not have a similar relation for the points $\Qunun$; we therefore
define the points $Q_{1}' , \ldots , Q_{n} '$ on the circle by requiring
that the analogue of (3.8) holds, i.e. by putting:
\begin{equation}
Q_{i} ' \ := \ T ( Q_{ \tau (i) } ), \ \ \ 1 \leq i \leq n.
\end{equation}
It is readily checked that $( \Punun ; Q_{1}' , \ldots , Q_{n} ')$ is an
$\ee '$-insertion bead, with $\ee '$ as in the statement of the lemma.

\newpage

\[
\mbox{{\bf Figure 2.}}
\]

With these notations, the inequality (3.6) that is to be proved has the
following geometric interpretation:
\begin{equation}
\left\{  \begin{array}{c}
{ \mbox{ If $[Q_{i},Q_{j}]$ intersects at least one $[P_{h},P_{k}]$ with
$h \ecpi k$, }  }  \\
                  \\
{ \mbox{ then $[ T (Q_{i}), T (Q_{j})]$ intersects at least one 
$[P_{h},P_{k}]$ with $h \ecpiprim k$. }  } 
\end{array}  \right.
\end{equation}
Indeed, the first line of (3.10) means $i \not\ecrho j$ (see 2.3), while 
the second line of (3.10) means $\tau (i) \not\ecrhoprim \tau (j)$
(where we also took (3.9) into account). Hence the statement in (3.10) is
that: $i \not\ecrho j \Rightarrow \tau (i) \not\ecrhoprim \tau (j)$
( $\Leftrightarrow i \stackrel{ \tau ( \rho ')}{\not\sim} j$ ), or in other 
words: $i \ectrp j \Rightarrow i \ecrho j;$ but this is exactly (3.6).

So we are left to prove (3.10). We shall divide the argument in 
three cases.

a) If $Q_{i}$ and $Q_{j}$ are on opposite sides of the line
$P_{t_{0}} P_{t_{1}}$. Then, as is clear from the definition 
of $T$, the points $T (Q_{i})$ and $T (Q_{j})$ are
still on opposite sides of $P_{t_{0}}P_{t_{1}}$; hence,
$[ T (Q_{i}) , T (Q_{j}) ]$ intersects
$[ P_{t_{0}} , P_{t_{1}} ],$ with $t_{0} \ecpiprim t_{1}.$

b) If $Q_{i}$ and $Q_{j}$ are on the same side of the line 
$X_{0}X_{1}.$ Consider $1 \leq h,k \leq n$ such that $h \ecpi k$ 
and $[Q_{i},Q_{j}] \cap [P_{h},P_{k}] \neq \emptyset$. From 
(3.7) it follows that $P_{h}$ and $P_{k}$ are also on the same 
side of $X_{0}X_{1}.$ But then, since $T$ is affine on any
of the two open half-planes determined by $X_{0}X_{1},$ we get:
\[
[ T (Q_{i}), T (Q_{j}) ] \cap 
[ P_{\tau (h)} , P_{\tau (k)} ] \ = \ 
T ( [ Q_{i} , Q_{j} ] ) \cap 
T ( [ P_{h} , P_{k} ] ) \ = \
T ( [ Q_{i} , Q_{j} ]  \cap [ P_{h} , P_{k} ] ) 
\ \neq \ \emptyset ;
\]
and in addition we have $\tau (h) \ecpiprim \tau (k)$, because
$h \ecpi k$ and $\pi ' = \tau ( \pi ).$

c) If $Q_{i}$ and $Q_{j}$ are on the same side of the line 
$P_{t_{0}}P_{t_{1}}$, but are on different sides of the line 
$X_{0}X_{1}.$ This is only possible if one of $i,j$ is in 
$\{ t_{0},t_{1} \}$ (say, for definiteness, that $i=t_{0}$), and 
the other one ($j$) is in $\nn \setminus \{ t_{0}, \ldots , t_{1} \}$. 
In order to have $Q_{i}$ and $Q_{j}$ on the same side of
$P_{t_{0}}P_{t_{1}}$ we must also assume that the $t_{0}$-th component
of $\ee$ is $1$; this implies that $Q_{i} = Y_{0},$ 
$T (Q_{i}) = Y_{1}.$

Consider now $1 \leq h,k \leq n$ such that $h \ecpi k$ and 
$[Q_{i},Q_{j}] \cap [P_{h},P_{k}] \neq \emptyset$. By using (3.7) 
we see that $h,k$ must necessarily be in $\nn \setminus \{ t_{0}, \ldots
,t_{1} \}$; note that we also have $h \ecpiprim k$, since $\pi ' =
\tau ( \pi )$ and $h,k$ are fixed by $\tau$. 

We know that $[ P_{h},P_{k} ] \cap [ Y_{0} , Q_{j} ] \neq \emptyset$
(previous paragraph), and that 
$[ P_{h},P_{k} ] \cap [ Y_{0} , Y_{1} ] = \emptyset$ (Eqn. (3.7)).
The line $P_{h}P_{k}$ cannot cross exactly one edge of the triangle
$Y_{0}Y_{1}Q_{j}$, therefore we must have 
$[ P_{h},P_{k} ] \cap [ Y_{1} , Q_{j} ] \neq  \emptyset$; but this means
exactly that 
$[ P_{h},P_{k} ] \cap [ T (Q_{i}) , T (Q_{j}) ] \neq  \emptyset$, 
and since $h \ecpiprim k$, the verification of (3.10) is complete.

At some points during the presentation of this proof, it was implicitly
assumed that $t_{0} \neq t_{1}$ and that $\nn \setminus \{ t_{0},
\ldots , t_{1} \} \neq \emptyset$; it is however easy to see that, with
the appropriate minor modifications, the argument is also working for
the two extreme cases when $t_{0}=t_{1}$ or $[t_{0},t_{1}] = [1,n].$
{\bf QED}

$\ $

{\bf 3.6 Lemma.} If $( \pi , \ee ) \in NCO(n) \times \strings$ and if
$\twtw ( \pi , \ee ) =: ( \pi ' , \ee ')$, then the numbers $d( \ee )$
and $d( \ee ')$ (counting how many 2's are in the $n$-tuples $\ee$ and
$\ee '$) have different parities.

$\ $

{\bf Proof.} The sequence $\ee '$ is obtained from $\ee$ by changing its 
components with indices in the twist interval $[t_{0},t_{1}]$ of $\pi$.
Thus if  $\ee \egdef ( \lunun )$ and if we put: $d'$ = 
$| \{ m \ | \ t_{0} \leq m \leq t_{1} , \ l_{m} = 2 \ \} | ,$ $d''$ = 
$| \{ m \ | \ m<t_{0}$ or $m>t_{1},$ and $l_{m} =2 \ \} |,$ then we obtain:
\[
d( \ee ) = d' + d'', \ \ \ d( \ee ') = ((t_{1}-t_{0}+1)-d') + d''.
\]
Hence $d( \ee ) + d( \ee ') = (t_{1}-t_{0})+1+2d''.$ Since $t_{0}$ and 
$t_{1}$ have the same parity (see Notations 3.4.1), we obtain that
$d( \ee ) + d( \ee ')$ is an odd number, and the conclusion follows.
{\bf QED}

$\ $

{\bf 3.7 Proposition.} Let $\ncps$ be a $\noncom$, and let $a,b \in \A$ be 
free with respect to $\varphi .$ Then for every $n \geq 1$ we have:
\begin{equation}
\sum_{ \begin{array}{c}
{\scriptstyle \pi \in NCO(n) }  \\
{\scriptstyle \ee \in \strings}
\end{array}  } \  (-1)^{d( \ee )} \
[ \mbox{coef } (n) ; \pi ] ( \rmua ) \cdot
[ \mbox{coef } (n) ; C_{\ee} ( \pi ) ] ( \mmub ) \ = \ 0.
\end{equation}

$\ $

{\bf Proof.} Let us denote the sum appearing in the left-hand side of (3.11)
by $L.$ By performing (in the named sum) the change of variable provided by
the bijection $\twtw$ of Definition 3.4.3, we obtain that:
\begin{equation}
L \ = \ \sum_{ \begin{array}{c}
{\scriptstyle \pi \in NCO(n) }    \\
{\scriptstyle \ee \in \strings} \\
{\scriptstyle \twtw ( \pi , \ee ) =: ( \pi ' , \ee ' )}
\end{array}  } \  (-1)^{d( \ee ' )} \
[ \mbox{coef } (n) ; \pi ' ] ( \rmua ) \cdot 
[ \mbox{coef } (n) ; C_{\ee ' } ( \pi ' ) ] ( \mmub ) .
\end{equation}
Hence we can write:
\begin{equation}
2L \ = \ \sum_{ \begin{array}{c}
{\scriptstyle \pi \in NCO(n) }    \\
{\scriptstyle \ee \in \strings} \\
{\scriptstyle \twtw ( \pi , \ee ) =: ( \pi ' , \ee ' )}
\end{array}  } \  \left\{ (-1)^{d( \ee ) }
[ \mbox{coef } (n) ; \pi ] ( \rmua ) \cdot 
[ \mbox{coef } (n) ; C_{\ee} ( \pi ) ] ( \mmub ) 
\right.
\end{equation}
\[
\left.
+ \  (-1)^{d( \ee ' )} \
[ \mbox{coef } (n) ; \pi ' ] ( \rmua ) \cdot 
[ \mbox{coef } (n) ; C_{\ee ' } ( \pi ' ) ] ( \mmub ) \ \right\} .
\]
But for every $( \pi , \ee ) \in NCO(n) \times \strings$, with
$\twtw ( \pi , \ee ) =: ( \pi ' , \ee ' ),$ it happens that
\begin{equation}
(-1)^{d( \ee ) } [ \mbox{coef } (n) ; \pi ] ( \rmua ) \cdot 
[ \mbox{coef } (n) ; C_{\ee} ( \pi ) ] ( \mmub )+\qquad\qquad\qquad
\qquad\  
\end{equation}
\[
+ \  (-1)^{d( \ee ' )} \ [ \mbox{coef } (n) ; \pi ' ] ( \rmua ) 
\cdot [ \mbox{coef } (n) ; C_{\ee ' } ( \pi ' ) ] ( \mmub ) \ = \ 0.
\]
Indeed, we have $[ \mbox{coef } (n) ; \pi ] ( \rmua )$ =
$[ \mbox{coef } (n) ; \pi ' ] ( \rmua )$ because $\pi$ and $\pi '$ have
the same block structure (obvious from (3.3));
we have $[ \mbox{coef } (n) ; C_{\ee} ( \pi ) ] ( \mmub )$ =
$[ \mbox{coef } (n) ;C_{\ee'}( \pi ') ]$ $ ( \mmub )$ because 
$C_{\ee} ( \pi )$ and $C_{\ee '} ( \pi ')$ have the same block structure 
(Lemma 3.5); and $(-1)^{d( \ee )} + (-1)^{d( \ee ')} =0,$
by Lemma 3.6.

>From (3.13) and (3.14) it is obvious that $L=0.$ {\bf QED}

$\ $

{\bf 3.8 Conclusion of the proof of Theorem 3.1.} Let $\ncps$ be a 
non-commutative probability space,
and let $a,b \in \A$ be free with respect to $\varphi$. By combining the 
Propositions 3.3 and 3.7 we obtain that, for every $n \geq 1:$
\begin{equation}
\varphi ( ( \comm )^{n} ) \ = \ 
\sum_{ \begin{array}{c}
{\scriptstyle \pi \in NCE(n) }  \\
{\scriptstyle \ee \in \strings}
\end{array}  } \  (-1)^{d( \ee )} \
[ \mbox{coef } (n) ; \pi ] ( \rmua ) \cdot 
[ \mbox{coef } (n) ; C_{\ee} ( \pi ) ] ( \mmub ).
\end{equation}
Note that on the right-hand side of (3.15) only the even free cumulants of
$a$ -- i.e. the coefficients of $\er ( \mu_{a} )$ -- are appearing, because 
one takes $[ \mbox{coef } (n); \pi ] ( \rmua )$ for 
$\pi$'s having only blocks of even size.

Consider now the situation discussed in Lemma 3.2, where we also have a
non-commutative probability space
$\pstild$ and $\atild , \btild \in \Atild,$ such that the hypothesis
(i)+(ii)+(iii$'$) are fulfilled. By writing the counterpart of Eqn. (3.15) for
$\atild$ and $\btild$, we clearly obtain that 
$\varphi ( ( \comm )^{n} )$ = $\varphi ( ( \comtild )^{n} )$
for every $n \geq 1,$ as desired. {\bf QED}

$\ $

$\ $

$\ $

\setcounter{section}{4}
\setcounter{equation}{0}
{\large\bf 4. $R$-diagonal pairs and free commutator} 

$\ $

If $a$ and $b$ are free in some non-commutative probability space, 
then $ab$ and $ba$ are of course not free; but if in addition we assume
that $a$ and $b$ are both even, then there still exists a definite and 
treatable relation between $ab$ and $ba$ -- namely that they form an
$R$-diagonal pair (notion introduced in \cite{NS3}). This fact has a
key role in our approach to the free commutator, and the main goal of
the present section is the presentation of its proof.
Before arriving to this (in Theorem 4.5), we shall review from \cite{NS3}
some facts about $R$-diagonality.

$\ $

{\bf 4.1 Definition.} 
Let $\ncps$ be a tracial \noncom and $x_1,x_2\in\A$.
We call $(x_1,x_2)$ an {\it $R$-diagonal pair} if its 
$R$-transform has the special form
\begin{equation}
\lb R(\mu_{x_1,x_2})\rb (z_1,z_2)=f(z_1z_2)+f(z_2z_1)
\end{equation}
for some series $f$ of only one variable; this $f$ is called the
{\it determining series} of the pair $(x_1,x_2)$.

$\ $

Two of the most important examples of $R$-diagonal pairs are $(c,c^*)$
and $(u,u^*)$ for $c$ a circular variable and $u$ a Haar unitary.
The determining series in these cases are given by $f=Id$ and
by $f=Moeb$, respectively (with $Id,$ $Moeb$ as in the comments following
the Definitions 2.6).

The main results of \cite{NS3} about $R$-diagonal pairs are
collected in the next theorem.

$\ $

{\bf 4.2 Theorem (\cite{NS3}).} 
Let $\ncps$ be a tracial \noncom. 

1) Let $(x_1,x_2)$ be an $R$-diagonal pair in $\ncps$. Then the 
determining series $f$ of $(x_1,x_2)$ is given by
\begin{equation}
f=R(\mu_{x_1x_2})\,\freestar\, Moeb.
\end{equation}

2) Let $x_1,x_2,p_1,p_2\in\A$ such that $(x_1,x_2)$ is an $R$-diagonal
pair and such that $\{x_1,x_2\}$ is free from $\{p_1,p_2\}$. Then
$(x_1p_1,p_2x_2)$ is also an $R$-diagonal pair.

3) Let, in the situation of part 2), $f$ and $g$ be the determining
series of $(x_1,x_2)$ and $(x_1p_1,p_2x_2)$, respectively.
Then these two series are related by
\begin{equation}
g=f \,\freestar\, R(\mu_{p_1p_2}).
\end{equation}

$\ $

Our interest in $R$-diagonal pairs came from the fact that some important
properties of the circular element have natural generalizations to
this framework. Two situations of this kind that we found are concerning
polar decompositions of certain $R$-diagonal pairs (see \cite{NS3},
Application 1.9), and the occurrence  of $R$-diagonal pairs as free
off-diagonal compressions (see \cite{N1}). We state the polar decomposition 
result precisely, because we will need to refer to it in Section 4.7.

$\ $

{\bf 4.3 Proposition (\cite{NS3}).}
Let $\ncps$ be a tracial probability space and assume that
$\A$ is a von Neumann algebra and $\ff$ a normal faithful state.
Consider $x\in\A$ such that $(x,x^*)$ is an $R$-diagonal pair and 
such that $\mbox{Ker }x=\{0\}$. Let $x=up$ be
the polar decomposition of $x$. Then it follows that $u$ is
a Haar unitary and that $\{u,u^*\}$ is free from $p$.

$\ $

The following characterization of
$R$-diagonal pairs will be used in the proof of Theorem 4.5.

$\ $

{\bf 4.4 Proposition.} 
Let $\ncps$ be a tracial \noncom and
$x_1,x_2\in\A$. Then the following two properties are equivalent:
\newline
a) The pair $(x_1,x_2)$ is $R$-diagonal.
\newline
b) If $u \in \A$ is a Haar unitary such that $\{u,u^*\}$ is free from 
$\{x_1,x_2\}$, then $\mu_{x_1,x_2}=\mu_{x_1u,u^*x_2}$.

$\ $

{\bf Proof.} $a)\Rightarrow b):$
By Theorem 4.2, the pair $(x_1u,u^*x_2)$ is $R$-diagonal and its
determining series $g$ is connected with the determining series $f$
of $(x_1,x_2)$ by
\begin{equation}
g=f\,\freestar\, 
R(\mu_{uu^*})=f\,\freestar\, R(\mu_1)=f\,\freestar\, Id \ = \ f.
\end{equation}
This implies that the $R$-transforms -- and hence also the 
distributions -- of the two pairs are the same.

$b)\Rightarrow a):$ The pair $(x_1u,u^*x_2)$ 
is $R$-diagonal by Theorem 4.2.2 and the fact that
$(u,u^*)$ is $R$-diagonal. Hence $(x_1,x_2)$ is also $R$-diagonal
(since it has the same distribution as $(x_1u,u^*x_2)$ ). {\bf QED}

$\ $

We now arrive to the result announced at the beginning of the section.

$\ $

{\bf 4.5 Theorem.}
Let $\ncps$ be a \noncom and consider $a,b\in\A$ such
that $a$ and $b$ are both even and such that $a$ is free from $b$.
Then $(ab,ba)$ is an $R$-diagonal pair.

$\ $

{\bf Proof.} Note first that we can replace $(\A,\ff)$ by $(\A',\ff')$,
where $\A'$ is the unital algebra generated by $a$ and $b$ and
$\ff'$ is the restriction of $\ff$ to $\A'$. Since $(\A',\ff')$
is always tracial (by the fact that the reduced free product of traces
gives again a trace, see e.g. Proposition 2.5.3 in \cite{VDN}), we
can assume without loss of generality that $(\A,\ff)$ is a tracial
non-commutative probability space.

Another assumption we can make (by embedding $\ncps$ into a larger
tracial non-commutative probability space) is that there exists a 
Haar unitary $u \in \A$ such that $\{ u, u^* \}$ is free from $\{ a,b \}$.

Put $x_1:=ab$, $x_2:=ba$, and $y_1:=x_1u=abu,$ $y_2:=u^*x_2=u^*ba.$
For a sequence $\varepsilon=(l_1,\dots,l_n)$ of 1's and 2's we will
use the abbreviations $x_\varepsilon:=x_{l_1}\dots x_{l_n}$,
$y_\varepsilon:=y_{l_1}\dots y_{l_n}$. By Proposition 4.4, what we have
to show is that $\mu_{x_1,x_2}=\mu_{y_1,y_2}$, or equivalently that: 
\begin{equation}
\varphi(x_\varepsilon) \ = \ \varphi(y_\varepsilon)
\end{equation}
for all sequences $\varepsilon$. For the rest of the proof we fix such 
a sequence $\varepsilon=(l_1,\dots,l_n)$ 
(with $n\in \N$ and $l_i\in \{1,2\}$), about which we will prove that 
(4.5) holds.

By Proposition 2.10, we have
\begin{equation}
\varphi(x_\varepsilon)=\sum_{\pi\in NC(n)}
\lb \mbox{coef }(n);\pi\rb(
R(\mu_a))\cdot \lb \mbox{coef }(n);C_\varepsilon(\pi)\rb
(M(\mu_b)).
\end{equation}
Since $a$ and $b$ are even, 
the odd cumulants of $a$ and the odd moments of $b$ vanish; therefore, 
by Proposition 2.4.2, only $\varepsilon$-alternating $\pi$'s give
a non-vanishing contribution to the sum.

On the other hand, by Proposition 2.13,
\begin{equation}
\varphi(y_\varepsilon)=\sum_{\pi\in NC(n)}
\lb \mbox{coef }(n);\pi\rb(R(\mu_a))\cdot
\lb \mbox{coef }
(\varepsilon);C_\varepsilon(\pi)\rb (M(\mu_{bu,u^*b})).
\end{equation}
The pair $(bu,u^*b)$ is $R$-diagonal (by Theorem 4.2), and in particular 
all the free cumulants of odd length of $\mu_{bu,u^*b}$ are equal to zero; 
as a consequence, the joint moments of odd length of $(bu,u^*b)$ are
all equal to zero, too (see Remark 2.12). But then Proposition 2.4.2 
implies again that we can restrict in our summation to $\ee$-alternating 
$\pi$'s.

So let us fix such an $\varepsilon$-alternating $\pi$; by comparing (4.6)
and (4.7), it only remains to show that
\begin{equation}
\lb \mbox{coef }(n);C_\varepsilon(\pi)\rb(M(\mu_b))=
\lb \mbox{coef }
(\varepsilon);C_\varepsilon(\pi)\rb (M(\mu_{bu,u^*b})).
\end{equation}
Finally, (4.8) is indeed true, because the partition $C_{\ee} ( \pi )$ is also
$\ee$-alternating (by Proposition 2.4.1) -- hence when we explicitly write
the right-hand side of (4.8) as a product (comp. Eqn. (2.14)), 
all the appearing $u$'s are cancelled by corresponding $u^*$'s.
{\bf QED}

$\ $

{\bf 4.6 Remarks.}

1) Both the $R$-diagonal pairs $(c,c^*)$ and $(u,u^*)$, with $c$ circular
and $u$ Haar unitary, can arise in the form $(ab,ba)$ with $a,b$ free and
both even; in fact both examples can be obtained when we ask in addition
that $b^2 =1,$ i.e. $\mu_b=1/2(\delta_{-1}+\delta_1)$. With this assumption,
the $R$-diagonal pair $(x_1 , x_2 ) := (ab,ba)$ has determining series
\begin{equation}
f=R(\mu_{abba})\,\freestar\, Moeb=R(\mu_{a^2})
\,\freestar\, Moeb .
\end{equation}
If we also ask that $a^2=1$, then $f=Moeb$ and $(x_1,x_2)$ is
a Haar unitary. If $a$ is a semicircular element, then $a^2$ is
a free Poisson element with $R(\mu_{a^2})=Zeta$; hence 
$f=Zeta \ \freestar \  Moeb = Id,$ and $(x_1,x_2)$ is a circular element.

2) The `dual version' of the statement of Theorem 4.5 is also true, that
is: if $\ncps$ is a tracial non-commutative probability space and if
$x_1,x_2\in\A$ are such that $(x_1,x_2)$ is an $R$-diagonal pair, then 
$x_1x_2$ and $x_2x_1$ are free. The proof of this assertion can be made 
as follows. Without loss of generality we can assume that there exists a
Haar unitary $u \in \A$ such that $\{u,u^*\}$ is free from
$\{x_1,x_2\}$. Then, by Proposition 4.4, $(x_1,x_2)$ has the
same distribution as $(x_1u,u^*x_2)$; so instead of proving that $x_1x_2$ 
is free from $x_2x_1$, it suffices to prove that $(x_1u)(u^*x_2)$ is free 
from $(u^*x_2)(x_1u).$ But the latter assertion is easily verified, by 
just using the definition of freeness.

\newpage

{\bf 4.7 Proofs of Propositions 1.10, 1.18.}

1) In 1.10 we have to show that if $a,b$ are free and both even in the
non-commutative probability space $\ncps$, then $\mu_{i(ab-ba)}$ = 
$\mu_{ab+ba}$. We saw in Theorem 4.5 that the pair $(ab,ba)$ is
$R$-diagonal, hence that
\[
[ R( \mu_{ab,ba} ) ] (z_{1},z_{2}) \ = \ f(z_{1}z_{2}) \ + \
f( z_{2}z_{1} ),
\]
with $f$ a series of one variable. But then the Eqns.(2.23-24) of 
Proposition 2.14 give us that
\[
\lb R(\mu_{i(ab-ba)}\rb(z) \ = \ 2f(z^2) \ = \
\lb R(\mu_{ab+ba})\rb(z),
\]
which entails (Remark 2.8) that $i(ab-ba)$ and $ab+ba$ have the same
distribution.

2) Proposition 1.18 is obtained by putting together the results in 4.3 and
4.5.  {\bf QED}

$\ $

$\ $

$\ $

\setcounter{section}{5}
\setcounter{equation}{0}
{\large\bf 5. Proof of the free commutator formula} 

$\ $

We shall first prove the Theorem 1.2 under the extra assumption that the
elements $a,b$ in its statement are even, and then in full 
generality -- in the Sections 5.3 and 5.4, respectively. Before going to 
the proof of 1.2, we need to put into evidence one more formula, concerning 
the $R$-transform of the square of an even element. This formula follows from
a natural bijection between $NCE(2n)$ and the set of intervals of the
poset $NC(n),$ which was observed in \cite{NS3}, Corollary 4.5:

$\ $

{\bf 5.1 Proposition (\cite{NS3}).} Recall that $NCE(2n)$
denotes the set of non-crossing partitions $\sigma$ of
$\{ 1,2, \ldots , 2n \}$, with the
property that every block of $\sigma$ has an even number of elements.
There exists a bijection $\Psi : NCE(2n) \ \rightarrow$
$\{ ( \pi , \rho ) \ |$ $\pi , \rho \in NC(n), \ \pi \leq \rho \}$,
such that if $\Psi ( \sigma ) = ( \pi , \rho ),$ then: $\sigma$ and $\pi$
have the same number of blocks, say $k;$ and moreover, one can write
$\sigma = \{ B_{1} , \ldots , B_{k} \}$,
$\pi = \{ A_{1} , \ldots , A_{k} \}$, in such a way that 
$|B_{j}| = 2|A_{j}|$ for every $1 \leq j \leq k.$
 
$\ $

{\bf 5.2 Proposition.} We have the formula
\begin{equation}
[ R( \mu_{a} ) ] (z) \ = \ [ \ R( \mu_{a^{2}} ) 
\,\freestar\, Moeb \ ] (z^{2}) ,
\end{equation}
or equivalently
\begin{equation}
\RE(\mu_a)=R(\mu_{a^2})\,\freestar\, Moeb,
\end{equation}
holding for every even element $a$ in some non-commutative probability
space $\ncps$.

$\ $

{\bf Proof.} Let us denote $[ R( \mu_{a} ) ] (z) \egdef \infsum \alpha_{n}
z^{n} .$ Due to the hypothesis that $a$ is even, we have $\alpha_{m} =0$
for every odd $m$ (Remark 2.8). Then the moment-cumulant formula (2.8)
gives us:
\begin{equation}
\varphi ( a^{2n} ) \ = \ \sum_{
\begin{array}{c} 
{\scriptstyle \sigma \in NCE(2n) }  \\
{\scriptstyle \sigma \egdef \{ B_{1} , \ldots , B_{k} \} }
\end{array}  }
\alpha_{|B_{1}|} \cdots \alpha_{|B_{k}|} , \ \ \  n \geq 1.
\end{equation}

We next use the bijection $\Psi$ of Proposition 5.1 to make a 
``change of
variable'' in the summation in (5.3), thus obtaining:
\[
\varphi ( a^{2n} ) \ = \ \sum_{
\begin{array}{c} 
{\scriptstyle \pi \leq \rho \,(\pi,\sigma\in NC(n)) }  \\
{\scriptstyle \pi \egdef \{ A_{1} , \ldots , A_{k} \} }
\end{array}  }
\alpha_{2|A_{1}|} \cdots \alpha_{2|A_{k}|} 
\]
\begin{equation}
= \ \sum_{
\begin{array}{c} 
{\scriptstyle \pi \in NC(n) }  \\
{\scriptstyle \pi \egdef \{ A_{1} , \ldots , A_{k} \} }
\end{array}  }
\alpha_{2|A_{1}|} \cdots \alpha_{2|A_{k}|} 
\cdot | \{ \rho \in NC(n) \ | \ \rho \geq \pi \} | .
\end{equation}

Since
\[
\RE(\mu_a)(z) \ = \ \infsum \alpha_{2n} z^{n} ,
\]
the product $\alpha_{2|A_{1}|} \cdots \alpha_{2|A_{k}|}$ appearing 
in (5.4) is just $[ \mbox{coef } (n) ; \pi ] (\RE(\mu_a))$.
The other factor showing up in the general term of the summation in (5.4),
$| \{ \rho \in NC(n) \ | \ \rho \geq \pi \} | ,$ can also be viewed as 
$| \{ \rho ' \in NC(n) \ | \ \rho ' \leq K( \pi ) \} | $ -- because the
Kreweras complementation map on $NC(n)$ sends 
$\{ \rho \in NC(n) \ | \ \rho \geq \pi \}$ bijectively onto
$\{ \rho ' \in NC(n) \ | \ \rho ' \leq K( \pi ) \}$. A direct inspection 
of how the operation $\freestar$ is defined (Equation (2.3)) gives us that 
the latter cardinality can be identified as $[ \mbox{ coef} (n); K( \pi )]
( Zeta \ \freestar \ Zeta )$.

Hence our formula for $\varphi ( a^{2n} )$ becomes:
\begin{equation}
\varphi ( a^{2n} ) \ = \ \sum_{\pi \in NC(n)}  
[ \mbox{coef } (n); \pi ] (\RE(\mu_a)) \cdot
[ \mbox{ coef} (n); K( \pi )] ( Zeta \ \freestar \ Zeta ), \ \ \ n \geq 1.
\end{equation}
But the right-hand side of (5.5) is exactly the coefficient of order $n$
in $\RE(\mu_a) \ \freestar \ Zeta \  
\freestar \ Zeta,$ again by the formula (2.3). 
Since $\varphi ( a^{2n} )$ can be, on the other hand, viewed as the 
coefficient of order $n$ in the moment series $M( \mu_{a^{2}} ),$ the 
conclusion we draw from this calculation is that:
\begin{equation}
M( \mu_{a^{2}} ) \ = \ 
\RE(\mu_a) \ \freestar \  Zeta \ \freestar  \ Zeta.
\end{equation}

Finally, we $\freestar$-operate with $Moeb \ \freestar \ Moeb$ on both
sides of (5.6). If we take into account that $Moeb$ is the inverse of $Zeta$
under $\freestar$, and that $M( \mu_{a^{2}} ) \ \freestar \ Moeb$ =
$R( \mu_{a^{2}} )$ (by Eqn. (2.6)), we see that (5.6) is turned into 
(5.2), which is equivalent to (5.1). {\bf QED}

$\ $

{\bf 5.3 Proof of Theorem 1.2 for even elements.} 
If, in the notations of 1.2, we assume 
in addition that $a$ and $b$ are even, then the pair $(ab,ba)$
is $R$-diagonal (Theorem 4.5); or in other words, the joint $R$-transform of
$ab$ and $ba$ has the form
\begin{equation}
[ R( \mu_{ab,ba} ) ] (z_{1} , z_{2} ) \ = \ h(z_{1}z_{2}) + h(z_{2}z_{1}),
\end{equation}
with $h$ the determining series of the pair $(ab,ba).$ From (5.7) and
and Eqn. (2.23) of 2.14 it follows that
\begin{equation}
[ R( \mu_{i(ab-ba)} ) ] (z) \ = \ h((iz)(-iz)) + h((-iz)(iz)) \ = \ 
2h( z^{2} );
\end{equation}
so what we actually have to prove is the equality:
\begin{equation}
h \ = \ \RE(\mu_a) \ \freestar \ \RE(\mu_b) \ \freestar \ Zeta.
\end{equation}
Now, Theorem 4.2 gives us the formula
\begin{equation}
h \ = \ R( \mu_{(ab)(ba)} ) \, \freestar\, Moeb.
\end{equation}
Without loss of generality we can assume that we are working in a tracial 
non-commutative probability space (compare beginning of proof of 
Theorem 4.5). Thus we have that $\mu_{abba} = \mu_{a^{2}b^{2}},$ and by
applying Equation (2.9) we get:
\begin{equation}
R( \mu_{abba} ) \ = \ R( \mu_{a^{2}b^{2}} ) \ = \ 
R( \mu_{a^{2}} ) \ \freestar \ R( \mu_{b^{2}} )
\end{equation}
(where we used, of course, that $a^{2}$ is free from $b^{2}$). 
Equations (5.10) and (5.11) imply together:
\[
h \ = \ R( \mu_{a^{2}} ) \ \freestar \ R( \mu_{b^{2}} ) \ \freestar \ Moeb
\]
\[
= \ ( R( \mu_{a^{2}} ) \ \freestar \ Moeb) \ \freestar \ Zeta \ \freestar \
( R( \mu_{b^{2}} ) \ \freestar \ Moeb ).
\]
\[
\stackrel{(5.2)}{ = } \ 
\RE(\mu_a)\,\freestar\, Zeta\,\freestar\, \RE(\mu_b) ,
\]
and thus (5.9) is obtained. {\bf QED}

$\ $

{\bf 5.4 Proof of Theorem 1.2 (general case).} Let $\ncps$ and
$a,b \in \A$ be as in Theorem 1.2. By a standard free product construction 
we can produce a second non-commutative probability space $\pstild$ and 
two elements $\atild , \btild \in \Atild$ which are free and both even, 
such that $\RE(\mu_a)=\RE(\mu_{\atild})$, 
$\RE(\mu_b)=\RE(\mu_{\btild})$. Theorem 3.1 gives us that 
$\mu_{\atild \btild - \btild \atild}$ (calculated in $\pstild$) coincides
with $\mu_{ab-ba}$ (calculated in $\ncps$). The result of the previous
section applies to $\atild$ and $\btild$, hence we just have to write: 
\[
\er ( \mu_{i(ab-ba)} ) \ = \
\er ( \mu_{i( \atild \btild - \btild \atild )} ) 
\ = \ 2( \RE(\mu_{\atild}) \ \freestar \ \RE(\mu_{\btild}) \ 
\freestar \ Zeta ) 
\]
\[
= \ 2( \RE(\mu_{a}) \ \freestar \ \RE(\mu_{b}) \ \freestar \ Zeta ). 
\mbox{\bf QED}
\]

$\ $

We now head towards the proofs of the corollaries of Theorem 1.2. As
mentioned in the Introduction, the proofs of 1.4 and 1.6 are obtained
by applying the combinatorial Fourier transform $\F$ in Equation (1.8)
of the theorem. Let us recall from \cite{NS1} that for a series 
$f(z) = \infsum \alpha_{n}z^{n}$ with $\alpha_{1} \neq 0,$ the series
$\F (f)$ is simply defined by
\begin{equation}
[ \F (f)](w) \ = \ \frac{1}{w} \ f^{<-1>} (w),
\end{equation}
where `$<-1>$' denotes inversion under composition. Recall also that
if $\mu : \C [X] \rightarrow \C$ is a distribution with $\mu (X) \neq 0,$
then its transform $S( \mu )$ is defined (\cite{V3}) as
\begin{equation}
[ S( \mu ) ](w) \ = \ \frac{1+w}{w} (M( \mu ))^{<-1>} (w),
\end{equation}
with $M( \mu )$ the moment series of $\mu$, as in Definition 2.7. The
combinatorial Fourier transform connects $R( \mu )$ with $S( \mu ),$ and
converts the $\freestar$-operation into multiplication, in the way shown
in the Eqns.(1.4--5) of the Introduction.

In the next proof we shall use the notation `$\circ D_{\lambda}$' for the
dilation of a formal power series with a $\lambda \in \C \setminus \{ 0 \}$;
more precisely, if $f(z) = \infsum \alpha_{n} z^{n},$ then
\begin{equation}
(f \circ D_{\lambda} ) (z) \ = \ f( \lambda z) \ := \
\infsum ( \alpha_{n} \lambda^{n} ) z^{n}.
\end{equation}
The $R$-transform and the $\freestar$-operation behave nicely with respect
to dilations, in the sense that we have the formulas:
\begin{equation}
R( \mu_{\lambda a} ) \ = \ R( \mu_{a} ) \circ D_{\lambda}
\end{equation}
for every $a$ in some $\ncps$ and for every $\lambda$ (see e.g. \cite{VDN},
Example 3.4.3); and
\begin{equation}
(f \circ D_{\lambda} ) \ \freestar \ g \ = \ 
f \ \freestar \ (g \circ D_{\lambda} ) \ = \ 
(f \ \freestar \ g) \circ D_{\lambda}   
\end{equation}
for every series $f,g$ and every $\lambda$ (see \cite{NS2}, Section 4.1).

$\ $

{\bf 5.5 Proof of Corollaries 1.4 and 1.6.} Let $\ncps$ be a 
non-commutative probability space, and let $a,b \in \A$ be free and
with non-zero variances $\gamma_{a} = \varphi (a^{2}) -
( \varphi (a))^{2},$ $\gamma_{b} = \varphi (b^{2}) -
( \varphi (b))^{2}.$ The fact that the variance of $c := i(ab-ba)$ 
satisfies $\gamma_{c} =2 \gamma_{a} \gamma_{b}$ comes out
by equating the coefficients of degree 1 on the two sides of Eqn. (1.8);
indeed, the linear coefficients of $\er ( \mu_{a} ), \er ( \mu_{b} ),
\er ( \mu_{c} )$ are exactly $\gamma_{a}, \gamma_{b}, \gamma_{c},$
respectively, and at the level of linear coefficients $\freestar$ is
just usual multiplication (see (2.3)).

It is convenient to rewrite (1.8) in the form:
\begin{equation}
\left\{  \begin{array}{l}
\er ( \mu_{c} ) \circ D_{1/ \gamma_{c}} \ = \ 2h \circ D_{1/2},  \\
                                                                 \\
h \ := \ ( \er ( \mu_{a} ) \circ D_{1/ \gamma_{a}} ) \ \freestar \
( \er ( \mu_{b} ) \circ D_{1/ \gamma_{b}} ) \ \freestar \ Zeta.
\end{array}  \right.
\end{equation}
By applying $\F$ in the two Equations (5.17), we obtain:
\[
[ \F ( \er ( \mu_{c} ) \circ D_{1/ \gamma_{c} } ) ] (w) \ = \ 
[ \F (h) ] ( \frac{w}{2} )
\]
(by just using the definition of $\F$ in (5.12)), and respectively
\[
\F (h) \ = \
\F ( \er ( \mu_{a} ) \circ D_{1/ \gamma_{a}} ) \cdot
\F ( \er ( \mu_{b} ) \circ D_{1/ \gamma_{b}} ) \cdot
\F ( Zeta ).
\]
If we also take into account that $[ \F ( Zeta )](w) = 1/(1+w),$
we thus arrive to:
\begin{equation}
[ \F ( \er ( \mu_{c} ) \circ D_{1/ \gamma_{c} } ) ] (w) \ = \ 
\frac{ 
[ \F ( \er ( \mu_{a} ) \circ D_{1/ \gamma_{a} } ) ] ( \frac{w}{2} ) \cdot 
[ \F ( \er ( \mu_{b} ) \circ D_{1/ \gamma_{b} } ) ] ( \frac{w}{2} )  }{
1+ \frac{w}{2}} .
\end{equation}

Now let us establish the formula:
\begin{equation}
[ \F ( \er ( \mu_{c} ) \circ D_{1/ \gamma_{c} } ) ] (w) \ = \ 
(1+w) \cdot [ S( \mu_{c^{2}/ \gamma_{c}} ) ] (w).
\end{equation}
Indeed, from Equation (5.2) we have, by also using (5.15--16),
\begin{equation}
\er ( \mu_{c} ) \circ D_{1/ \gamma_{c}}  \ = \
( R( \mu_{a^{2}} ) \circ D_{1/ \gamma_{c}} ) \ \freestar \ Moeb \ = \
R( \mu_{a^{2}/ \gamma_{c}} ) \ \freestar \ Moeb;
\end{equation}
then by applying $\F$ in (5.20) and by substituting
$\F( R( \mu_{a^{2} / \gamma_{c}} ))$ =
$S( \mu_{a^{2} / \gamma_{c}} )),$ $[ \F ( Moeb )](w) = 1+w,$ we obtain
(5.19).

By putting (5.18) and (5.19) together, we obtain:
\begin{equation}
[ S( \mu_{c^{2}/ \gamma_{c}} ) ] (w) \ = \
\frac{ 
[ \F ( \er ( \mu_{a} ) \circ D_{1/ \gamma_{a} } ) ] ( \frac{w}{2} ) \cdot 
[ \F ( \er ( \mu_{b} ) \circ D_{1/ \gamma_{b} } ) ] ( \frac{w}{2} ) }{ 
(1+ \frac{w}{2})(1+w) } .
\end{equation}

If $a$ and $b$ happen to be even, then we can write the counterparts of
(5.19) for $a$ and $b$, and substitute them into (5.21);
this leads exactly to the Equation (1.25) of Corollary 1.6. 

If $a$ and $b$ are not assumed to be even, then we just replace the $\F$'s 
on the right-hand side of (5.21) from their definition (Eqn. (5.12)); then
(5.21) becomes:
\begin{equation}
[ S( \mu_{c^{2} / \gamma_{c}} )] (w) \ = \ 
\frac{4 \gamma_{a} \gamma_{b}}{w^{2} (1 + \frac{w}{2})(1+w)} \cdot 
[ \er ( \mu_{a} ) ]^{<-1>} ( \frac{w}{2} ) \cdot
[ \er ( \mu_{b} ) ]^{<-1>} ( \frac{w}{2} ) .
\end{equation}
In (5.22) we multiply both sides with $w/(1+w),$ and then take their inverse
under composition; in this way the left-hand side of the equation becomes
$M( \mu_{ c^{2} / \gamma_{c} } )$ (see (5.13)). More precisely, we get the 
formula:
\begin{equation}
[ M( \mu_{c^{2} / \gamma_{c}} )] (z) \ = \ 
\left( \frac{4 \gamma_{a} \gamma_{b}}{w (1 + \frac{w}{2})(1+w)^{2} } \cdot 
[ \er ( \mu_{a} ) ]^{<-1>} ( \frac{w}{2} ) \cdot
[ \er ( \mu_{b} ) ]^{<-1>} ( \frac{w}{2} ) \right)^{<-1>} (z).
\end{equation}
If, finally, we do a dilation with $\gamma_{c} = 2 \gamma_{a} \gamma_{b}$ in
(5.23), and replace $z$ by $z^{2}$, then Equation (1.11) of Corollary 1.4 is
obtained. {\bf QED}

$\ $

{\bf 5.6 Proof of Corollary 1.13.} Let $a,b$ be free selfadjoint elements
in a $C^{*}$-probability space $\ncps$, such that $\mu_{a} = \nu_{1}$ and 
$\mu_{b} = \nu_{2}.$ Without loss of generality we can assume that $\varphi$
is a trace and (by appropriately enlarging $\ncps$) that there exists 
$d=d^{*} \in \A$ with $\mu_{d} = \frac{1}{2} ( \delta_{1} + \delta_{-1})$
and such that $d$ is classically independent (not free!) from $\{ a,b \}$.

We know, by (5.8),(5.10), that
\begin{equation}
[ R(\miab) ] (z)=h(z^2) \qquad \mbox{with} \qquad
h=R( \mu_{abba} )\,\freestar\, Moeb.
\end{equation}
We define $\widehat a\egdef \sqrt{abba} \cdot d$. Then $\widehat a$ is 
selfadjoint and even, and has $\widehat a^2=abba$. Equations (5.1) and 
(5.24) imply together that
\[
\lb R(\mu_{\widehat a})\rb(z)=\lb R( \mu_{abba} )\,\freestar\, Moeb\rb 
(z^2)=h(z^2) =[ R(\miab) ] (z),
\]
hence that
\begin{equation}
\lb\nu_1\bbox\nu_2\rb^{\bplus 1/2}=\miab=\mu_{\widehat a}.
\end{equation}
Since $\widehat a$ is a selfadjoint element in a $C^{*}$-probability space, 
its distribution is a probability measure; thus (5.25) implies part 1) of 
Corollary 1.13. Moreover, formula (1.29) in the part 2) of 1.13
also follows from (5.25), via the calculation:
\[
R(Q(\lb\nu_1\bbox\nu_2\rb^{\bplus 1/2})) \ = \ 
R(Q(\mu_{\widehat a})) \ = \
R(\mu_{\widehat a^2})
\]
\[
=R(\mu_{abba}) =R(\mu_{a^2}\,\btimes\,\mu_{b^2})=R(Q(\mu_a)\,\btimes\,
Q(\mu_b)).  \mbox{\bf QED}
\]

$\ $

We are  only left to present the proof of Corollary 1.16. In order to do
this, we introduce one more notation concerning commutator expressions:

$\ $

{\bf 5.7 Notation.} To each commutator expression $f$ of $n$ arguments
(in the sense of 1.15) we assign a {\em $\bbox$-depth vector} 
$(t_1,\dots, t_{n-1})$, in the following recursive way.
\newline
i) $f(\nu_1)=\nu_1$ has $\bbox$-depth $\emptyset$
\newline
ii) If $f(\nu_1,\dots,\nu_n)=\lb f_1(\nu_1,\dots,\nu_k)
\bbox f_2(\nu_{k+1},\dots,\nu_n)\rb$ and if
$f_1$ has $\bbox$-depth $(t_1,\dots,t_{k-1})$ and $f_2$
has $\bbox$-depth $(t_{k+1},\dots,t_{n-1})$, then $f$ has
$\bbox$-depth $(t_1+1,\dots,t_{k-1}+1,1,t_{k+1}+1,\dots,t_{n-1}+1)$.

$\ $

The $\bbox$-depth vector records how many brackets we have to cross in
order to reach the various `$\bbox$' signs inside the commutator expression
(in the same way as the depth vector of 1.15.2 was doing this for the 
arguments of the expression). We could afford not to mention the $\bbox$-depth
in the statement of Corollary 1.16, due to the following fact.

$\ $

{\bf 5.8 Lemma.} Let $f$ and $\widehat f$ be two commutator expressions
of $n$ arguments. If the depth vector of $f$ differs from the one of 
$\widehat f$ only by a permutation, then the same is true for the 
$\bbox$-depth vectors of $f$ and $\widehat{f}$.

$\ $

{\bf Proof.} By induction on $n$. For $n=1$ and $n=2$ the 
assertion is trivial. Let us assume we have proved it for $n-1$ and 
consider the statement for $n$. We denote the depth and $\bbox$-depth of $f$ 
by $(d_{1}, \ldots , d_{n})$ and $(t_{1}, \ldots,t_{n-1})$, and those of
$\widehat{f}$ by $( \widehat{d}_{1}, \ldots , \widehat{d}_{n})$ and
$( \widehat{t}_{1}, \ldots , \widehat{t}_{n-1} ),$ respectively.

Let $m$ be the biggest depth appearing in $f$ and $\widehat f$, i.e.
$$m:=\max\{d_1,\dots,d_n\}=\max\{\widehat d_1,\dots,\widehat d_n\}.$$
It is clear that there exist $k$ and $l$ such that
$d_k=d_{k+1}=m=\widehat d_l=\widehat d_{l+1}$ and that $t_k=m=\widehat t_l$.
Furthermore, there are commutator expressions $g$ and $\widehat g$ 
of $n-1$ arguments such that $f$ and $\widehat f$ are of the form
$$f(\nu_1,\dots,\nu_n)=g(\nu_1,\dots,\nu_{k-1},\lb\nu_k\bbox\nu_{k+1}\rb,
\nu_{k+2},\dots,\nu_n)$$
and
$$\widehat f(\nu_1,\dots,\nu_n)
=\widehat g(\nu_1,\dots,\nu_{l-1},\lb\nu_l\bbox\nu_{l+1}\rb,
\nu_{l+2},\dots,\nu_n) . $$
The depths of $g$ and $\widehat g$ are determined by the depths
of $f$ and $\widehat f$, respectively, just by replacing $(\dots,m,m,\dots)$
by $(\dots,m-1,\dots)$. Thus the depth vectors of $g$ and $\widehat g$ differ
also only by a permutation; this implies, by the induction hypothesis, that
the same is true for the $\bbox$-depth vectors of $g$ and $\widehat g$.
But the $\bbox$-depths of $f$ and $\widehat f$ differ from those of $g$
and $\widehat g$ just by $t_k$ and $\widehat t_l$, respectively.
Since $t_k=\widehat t_l$, the assertion for $f$ and $\widehat f$ follows.
{\bf QED}

$\ $

{\bf 5.9 Proof of Corollary 1.16.} For $f,g$ formal power series without
constant coefficient (in 1 variable), and for $\lambda \in \C \setminus
\{ 0 \}$, we have the formula:
\begin{equation}
( \lambda f) \ \freestar \ ( \lambda g) \ = \ \lambda (f \ \freestar \ g)
\circ D_{\lambda},
\end{equation} 
where `$\circ D_{\lambda}$' is as defined in Eqn. (5.14) (see \cite{NS2},
Lemma 4.4). By using (5.26) twice, we can rewrite the Equation (1.8) of 
Theorem 1.2 in the alternative form:
\begin{equation}
\er ( \nu_{1} \bbox \nu_{2} ) \ = \ 
(2 \er ( \nu_{1} )) \  \freestar \ (2 \er ( \nu_{2} )) \  \freestar \ 
(2 Zeta) \circ D_{1/4},
\end{equation}
with $\nu_{1}, \nu_{2}$ arbitrary distributions. The recursive use of (5.27)
is very convenient for calculating the $\er$-transforms of higher order free
commutator expressions. More precisely, as is easily checked by induction, 
we obtain the following general formula: let $f$ be a commutator expression 
of $n$ arguments, with depth vector $( d_{1}, \ldots , d_{n} )$ and 
$\bbox$-depth vector $(t_{1}, \ldots , t_{n-1})$; if 
$\nu_{1}, \ldots, \nu_{n}$ are arbitrary
distributions, and if we denote $\nu := f( \nu_{1}, \ldots , \nu_{n}),$ then:
\begin{equation} 
\ER( \nu )=\lb(2^{d_1}\ER(\nu_1))\,\freestar\,\dots\,\freestar
(2^{d_n}\ER(\nu_n))\,\freestar\, (2^{t_1} Zeta)\,\freestar\,
\dots\,\freestar\, (2^{t_{n-1}} Zeta)\rb \circ D_{4^{-(n-1)}}.
\end{equation}
But, if we also take the Lemma 5.8 into account, the assertion of 
Corollary 1.16 is immediately implied by Eqn. (5.28). {\bf QED}

$\ $

$\ $

{\large \bf Acknowledgments.}
We are grateful to Dan Voiculescu for directing our attention to the 
problem of the free commutator, and to Philippe Biane and Phil Hanlon for 
useful discussions. The paper was completed while the first-named author
was visiting the Forschungsgruppe of Joachim Cuntz at Heidelberg, and
he would like to acknowledge the pleasant atmosphere of this visit.

$\ $

$\ $


\begin{thebibliography}{99}


\bibitem{Akh} N.I. Akhiezer.
The classical moment problem and some related questions in
analysis, Oliver and Boyd, 1965.

\bibitem{And} G.E. Andrews.
The theory of partitions, Encyclopedia of Mathematics and Its
Applications, Volume 2, Addison-Wesley, 1976.

\bibitem{BV} H. Bercovici, D. Voiculescu.
Free convolution of measures with unbounded support,
Indiana Univ. Math. J. 42 (1993), 733-773.

\bibitem{BV2} H. Bercovici, D. Voiculescu.
Regularity questions for free convolution, preprint 1996.

\bibitem{Bia1} Ph. Biane. 
Processes with free increments, Math. Z. (to appear).

\bibitem{Bia2} Ph. Biane.
The density of free stable distributions, preprint 1996.

\bibitem{K} G. Kreweras.
Sur les partitions non-croisees d'un cycle, Discrete Math.
1 (1972), 333-350.

\bibitem{N} A. Nica.
$R$-transforms of free joint distributions and non-crossing 
partitions, J. Funct. Anal. 135 (1996), 271-296. 

\bibitem{N1} A. Nica.
$R$-diagonal pairs arising as free off-diagonal compressions, 
Indiana Univ. Math. J. (to appear).

\bibitem{NS1} A. Nica, R. Speicher.
A ``Fourier transform'' for multiplicative functions on non-crossing
partitions, J. Algebraic Combinatorics (to appear).

\bibitem{NS2} A. Nica, R. Speicher.
On the multiplication of free $n$-tuples of non-commutative random
variables, with an Appendix by Dan Voiculescu: Alternative proofs for
the type II free Poisson variables and for the free compression results.
Amer. J. Math. 118 (1996), 799-837.

\bibitem{NS3} A. Nica, R. Speicher.
$R$-diagonal pairs -- a common approach to Haar unitaries and circular
elements, funct-an 9604012, Fields Institute Communications (to appear).

\bibitem{S1} R. Speicher.
Multiplicative functions on the lattice of non-crossing partitions
and free convolution, Math. Annalen 298 (1994), 611-628.

\bibitem{V1} D. Voiculescu.
Symmetries of some reduced free product C*-algebras, in
Operator algebras and their connection with topology and ergodic
theory, H. Araki et al. editors
(Springer Lecture Notes in Mathematics, volume 1132, 1985), 556-588.

\bibitem{V2} D. Voiculescu.
Addition of certain non-commuting random variables, J. Functional
Analysis 66 (1986), 323-346.

\bibitem{V3} D. Voiculescu.
Multiplication of certain non-commuting random variables,
J. Operator Theory 18 (1987), 223-235.

\bibitem{ICM} D. Voiculescu.
Free probability theory: random matrices and von Neumann algebras,
in Proceedings of the ICM 1994 (Birkh\"auser, 1995), 227-241.

\bibitem{VDN} D. Voiculescu, K. Dykema, A. Nica.
Free random variables, CRM Monograph Series, volume 1, AMS, 1992.

\end{thebibliography}
\end{document}